\documentclass[reqno]{amsart}
\usepackage[T1]{fontenc}
\usepackage{dsfont}
\usepackage{mathrsfs}
\usepackage[colorlinks, citecolor=blue, linkcolor=blue]{hyperref}
\usepackage{xcolor}
\usepackage[a4paper,asymmetric]{geometry}
\usepackage{mathscinet}
\usepackage{latexsym}
\usepackage{amsthm}
\usepackage{amssymb}
\usepackage{amsfonts}
\usepackage{amsmath}
\usepackage{longtable}
\usepackage{graphicx}
\usepackage{multirow}
\usepackage{multicol}
\usepackage{amsfonts, amsmath}
\usepackage{latexsym,bm,amsfonts,amssymb,pifont,mathbbol,bbm}
\usepackage{verbatim}

\newtheorem{theorem}{Theorem}[section]
\newtheorem{thm}[theorem]{Theorem}
\newtheorem{lemma}[theorem]{Lemma}
\newtheorem{lem}[theorem]{Lemma}
\newtheorem{remark}[theorem]{Remark}
\newtheorem{proposition}[theorem]{Proposition}
\newtheorem{prop}[theorem]{Proposition}
\newtheorem{corollary}[theorem]{Corollary}
\newtheorem{hyp}[theorem]{HYPOTHESIS}
\theoremstyle{definition}

\newtheorem{defn}[theorem]{Definition}

 \newtheorem{nott}{Notation}
\theoremstyle{remark}
\numberwithin{equation}{section}

 \DeclareMathAlphabet{\mathpzc}{OT1}{pzc}{m}{it}
 \DeclareMathAlphabet{\mathsfsl}{OT1}{cmss}{m}{sl}

  \newcommand{\FH}{\mathfrak{H}}

\newcommand{\dif}{\mathrm{d}}

\newcommand{\abs}[1]{\left\vert#1\right\vert}
\newcommand{\set}[1]{\left\{#1\right\}}

\newcommand{\norm}[1]{\left\Vert#1\right\Vert}

\newcommand{\E}{\mathbb{E}}

 \newcommand{\Rnum}{\mathbb{R}}

 \newcommand{\innp}[1]{\langle {#1}\rangle}
\newcommand{\Be}{\begin{equation}}
\newcommand{\Ee}{\end{equation}}
\newcommand{\Bs}{\begin{split}}
\newcommand{\Es}{\end{split}}
\newcommand{\Bes}{\begin{equation*}}
\newcommand{\Ees}{\end{equation*}}
\newcommand{\BT}{\begin{thm}}
\newcommand{\ET}{\end{thm}}
\newcommand{\Bp}{\begin{proof}}
\newcommand{\Ep}{\end{proof}}
\newcommand{\BL}{\begin{lem}}
\newcommand{\EL}{\end{lem}}
\newcommand{\BP}{\begin{proposition}}
\newcommand{\EP}{\end{proposition}}
\newcommand{\BC}{\begin{corollary}}
\newcommand{\EC}{\end{corollary}}
\newcommand{\BR}{\begin{remark}}
\newcommand{\ER}{\end{remark}}
\newcommand{\BD}{\begin{defn}}
\newcommand{\ED}{\end{defn}}
\newcommand{\BI}{\begin{itemize}}
\newcommand{\EI}{\end{itemize}}

\allowdisplaybreaks

\begin{document}
\title[Parameter estimation for Ornstein-Uhlenbeck processes]{Parameter estimation for an Ornstein-Uhlenbeck Process driven by a general Gaussian noise with Hurst Parameter $H\in (0,\frac12)$}
\author[Y. Chen]{Yong CHEN}
 \address{College of Mathematics and Information Science, Jiangxi Normal University, Nanchang, 330022, Jiangxi, China}
 \author[X. Gu]{Xiangmeng GU}
 \address{School of Mathematics and Statistics, Jiangxi Normal University, Nanchang, 330022, Jiangxi, China}
  \author[Y. Li]{Ying LI}
 \address{School of Mathematics and Computional Science, Xiangtan University, Xiangtan, 411105, Hunan, China. (Corresponding author.)}  \email{liying@xtu.edu.cn}
\begin{abstract}
In Chen and Zhou 2021, they consider an inference problem for an Ornstein-Uhlenbeck process driven by a general one-dimensional centered Gaussian process $(G_t)_{t\ge 0}$. The second order mixed partial derivative of the covariance function $ R(t,\, s)=\mathbb{E}[G_t G_s]$ can be decomposed into two parts, one of which coincides with that of fractional Brownian motion and the other is bounded by $(ts)^{H-1}$ with $H\in (\frac12,\,1)$, up to a constant factor.
 In this paper, we investigate the same problem but with the assumption of $H\in (0,\,\frac12)$. It is well known that there is a significant difference  between the Hilbert space associated with the fractional Gaussian processes in the case of $H\in (\frac12, 1)$ and that of $H\in (0, \frac12)$. The starting point of this paper is a new relationship between the  inner product of $\mathfrak{H}$ associated with the Gaussian process $(G_t)_{t\ge 0}$ and that of the Hilbert space $\mathfrak{H}_1$ associated with the fractional Brownian motion $(B^{H}_t)_{t\ge 0}$. Then we prove the strong consistency with $H\in (0, \frac12)$, and the asymptotic normality and the Berry-Ess\'{e}en bounds with $H\in (0,\frac38)$ for both the least squares estimator and the moment estimator of the drift parameter constructed from the continuous observations.  A good many inequality estimates are involved in and we also make use of the estimation of the inner product based on the results of $\mathfrak{H}_1$ in Hu, Nualart and Zhou 2019. 



\end{abstract}
\maketitle

\section{ Introduction}\label{sec 03}

In \cite{CZ 21},  the statistical inference problem of the unknown parameter $\theta$  is considered for the Ornstein-Uhlenbeck process defined by the following stochastic differential equation (SDE)
\begin{equation}\label{fOU}
\mathrm{d} X_t= -\theta X_t\mathrm{d} t+ \sigma \mathrm{d}G_t,\quad  t \in [0,T], \ T >0
\end{equation} where $X_0=0$ and $(G_t)_{t\ge 0}$ is a general one-dimensional centered Gaussian process satisfying the following Hypothesis \ref{hypthe 0001chenzhou}.  
\begin{hyp}\label{hypthe 0001chenzhou}
For $H\in (\frac12, \,1)$, the covariance function $ R(t,\, s)=\mathbb{E}[G_t G_s]$ for any $t\neq s \in [0,\infty)$ satisfies
\begin{align}\label{hyp 0001 first}
\frac{\partial^2}{\partial t\partial s}R(t,s)=  {H}(2H-1)\abs{t-s}^{2H -2}+\Psi(t,\,s),
\end{align}with 
\begin{align}\label{hyp 0001 second}
\abs{\Psi(t,\,s)}\le C'_{H} \abs{ts}^{H -1},
\end{align} where the constant $ C'_{H}\ge 0$ do not depend on $T$.
Moreover, for any $t \geq 0$, $R(0,t)=0$.
\end{hyp}
Without loss of generality,  $\sigma=1$ is also assumed. Suppose that only one trajectory $(X_t, t \in [0,T])$ can be obtained. When $\theta>0$,  i.e., in the ergodic case, the least squares estimator (LSE) and the moment estimator (ME) are respectively constructed from the continuous observations as follows:
 \begin{align}
\hat{\theta}_T&=-\frac{\int_0^T X_t\mathrm{d}X_t}{\int_0^T X_t^2\mathrm{d} t}=\theta-\frac{\int_0^T X_t\mathrm{d}G_t}{\int_0^T X_t^2\mathrm{d} t}, \label{hattheta}\\
\tilde{\theta}_{T}&=\Big( \frac{1}{ {H} \Gamma(2H ) T} \int_0^T X_t^2\mathrm{d} t \Big)^{-\frac{1}{2H}},\label{theta tilde formula}
\end{align}where the integral with respect to $G$ is interpreted in the Skorohod sense (or say a divergence-type integral).
 Then the strong consistency, the asymptotic normality,  and the Berry-Ess\'{e}en bounds are obtained in \cite{CZ 21}.  

In fact, the statistical inference problem about the parameter $\theta$ has been intensively studied over the past decades (see \cite{k04}, \cite{ls01} and the references therein) when the Gaussian process is Brownian motion. In the fractional Brownian motion case, the consistency property for the maximum likelihood estimation (MLE) method was obtained in \cite{KLB 02}, \cite{tv 07}, and the central limit theorem was proved in \cite{bcs11}, \cite{bk10}, and the LSE and the ME and their asymptotic behavior were studied in \cite{hu Nua 10}, \cite{hu nua zhou 19}.  The ME in the case of general stationary-increment Gaussian processes was considered in \cite{SV 18}. The MLE in the case of sub-fractional Brownian motion case was investigated in \cite{dmm11} and, recently,  the LSE in the case of mixed sub-fractional Brownian
motion was studied in \cite{cai xiao 18}.  We would like to mention some work for the non-ergodic case as well, i.e., $\theta<0$.  MLE was studied in \cite{bs83}, \cite{dk03} and the limiting distribution is Cauchy for the Brownian motion case, and the LSE in the case of fractional Brownian motion and other Gaussian processes was considered in \cite{beo11}, \cite{eeo16}, \cite{mendy 13}, \cite{AAE 20} and the references therein.  We also mention some work for the Ornstein-Uhlenbeck process driven by the non-Gaussian Hermite processes with periodic mean in \cite{ST 20}, \cite{Shen Yu Tang 21}  and the references therein.

Clearly, the following four types of fractional Gaussian processes satisfy Hypothesis ~\ref{hypthe 0001chenzhou}. 
The covariance function of the fractional Brownian motion $\{B^H_t, t\geq 0\}$ is
\begin{equation}\label{exmp1}
R^{B}(s,t)= \frac{1}{2}(s^{2H} + t^{2H} - |t-s|^{2H}).
\end{equation}
 The sub-fractional Brownian motion $\{S^H_t, t \geq 0\}$ with parameter $H\in (0,1)$ has the covariance function
$$R(t,s)=s^{2H}+t^{2H}-\frac{1}{2}\left((s+t)^{2H}+|t-s|^{2H}\right).$$
 The bi-fractional Brownian motion $\{B^{H,K}_t, t\geq 0\}$ with parameters $H,K \in (0, 1)$ has the covariance function
$$R(t,s)=\frac{1}{2^K}\left((s^{2H}+t^{2H})^K - |t-s|^{2HK}\right) .$$
The generalized sub-fractional Brownian motion $\set{S^{H,K}_t, t\ge 0} $ with parameters $H \in (0, 1),\,K \in[1,2)$ and $HK\in (0,1)$  has the covariance function 
$$ R(t,\, s)= (s^{2H}+t^{2H})^{K}-\frac12 \big[(t+s)^{2HK} + \abs{t-s}^{2HK} \big].$$
The Hurst parameters should be understood as $HK$ for both of the last two Gaussian processes.

When $\theta>0$, i.e., in the ergodic case, for the above LSE \eqref{hattheta} and ME \eqref{theta tilde formula}, most of the results in literatures are restricted to the Hurst parameter $H\in (\frac12, 1)$ except \cite{hu nua zhou 19} and \cite{CL 21} as far as we know. This is mainly due to the significant difference between the Hilbert space associated with the fractional Gaussian processes and the representation of their inner products in the case of $H\in (\frac12, 1)$ and that in the case of $H\in (0, \frac12)$. We would like to point out a remarkable fact that the monotonicity of the norm may not hold in the case of $H\in (0,\,\frac12)$, please see Remark~\ref{rem dandiao} in Section~\ref{sec pre} or refer to \cite{Jolis2007} for details.

However, it is clearly that for all the above four types fractional Gaussian processes, the identity \eqref{hypthe 0001chenzhou} and the inequality \eqref{hyp 0001 second} are valid for both $H\in  (0,\frac12)$ and $H\in  (\frac12,\,1)$. Then the question naturally arises if the asymptotic properties of  LSE \eqref{hattheta} and ME \eqref{theta tilde formula} are also valid for $H\in  (0,\frac12)$.  We have pointed out in the last paragraph that this problem have been solved in \cite{hu nua zhou 19} and \cite{CL 21}  in the case of fractional Brownian motion. In this paper, we will partly give an affirmative answer to this question in the case of general fractional Gaussian processes. For simplicity, we would like to discuss only the ergodic case, i.e., $\theta>0$. 

The starting point of the present paper is to establish the key inequalities \eqref{inequality 29} and \eqref{inner product 00.ineq} which relate the inner product of the Hilbert space associated with the general  fractional Gaussian process to that associated with the fractional Brownian motion. To this aim, we need improve Hypothesis ~\ref{hypthe 0001chenzhou} into the following form \cite{CDL 21}:
 \begin{hyp}\label{hypthe 1}
For $H \in  (0,\frac12)\cup (\frac12,1)$,  the covariance function $R(t,s)=\mathbb{E}[G_{t}G_{s}]$ satisfies that
 \begin{enumerate}
	\item  for any $s\geq 0$, $R(0,s)=0$.
	\item  for any fixed $s\in (0,T)$, $ R(t,s) $ is a continuous function on $[0,T]$ which is differentiable with respect to $t$ in $(0,s)\cup(s,T)$, such that $\frac{\partial  }{\partial t }R(t,s) $ is absolutely integrable.
	\item  
for any fixed $t\in (0,T)$, the difference
$$\frac{\partial R(t,s)}{\partial t} - \frac{\partial R^{B}(t,s)}{\partial t} $$
is a continuous function on $[0,T]$ which is differentiable with respect to $s$ in $(0,T)$ such that $\Psi (t,s)$, the partial derivative with respect to $s$ of the difference, satisfies \begin{align}\label{cond hyp2}
    |\Psi (t,s)|&\leq C_{H }^{^{\prime }}|ts|^{H -1},
\end{align}
where   
the constant $C_{H }^{^{\prime }} \ge 0$ do not dependent on $T$, and $R^B(t,s)$ is the covariance function of the fractional Brownian motion given in \eqref{exmp1}.
	 \end{enumerate}
\end{hyp}
It is easy to see that fractional Brownian motion, sub-fractional Brownian motion, bi-fractional Brownian motion, generalized sub-fractional Brownian motion and some other Gaussian processes are special examples to satisfy the Hypothesis \ref{hypthe 1}. For example,  the mixed Gaussian process   \cite{cai xiao 18} which is a linear combination of independent centered Gaussian processes satisfies Hypothesis \ref{hypthe 1} as long as each Gaussian process satisfies it. In this case, the mixed Gaussian process fails to be self-similar.

In this paper, we will prove the strong consistency and the central limit theorems for the two estimators. The Berry-Ess\'{e}en bounds will be also obtained. These results are stated in the following theorems.

\begin{thm}\label{thm strong}
Let the least squares estimator $\hat{\theta}$ and the moment estimator  $\tilde{\theta}_{T}$ be as \eqref{hattheta} and \eqref{theta tilde formula} respectively.
When Hypothesis~\ref{hypthe 1} is satisfied, both  $\hat{\theta}$ and  $\tilde{\theta}_{T}$ are strongly consistent, i.e.,
\begin{align*}
\lim_{T\to\infty}\hat{\theta}_T=\theta,\qquad \lim_{T\to\infty}\tilde{\theta}_T=\theta, \qquad a.s..
\end{align*}
\end{thm}
A famous theorem of Pickands, see e.g. \cite{Pik 69},  states that if the covariance $\mathrm{Cov}(\xi_0,\,\xi_t)$ of a stationary Gaussian process $\xi_t$ with unit variance satisfies 
\begin{align*}
\mathrm{Cov}(\xi_0,\,\xi_t)=1-c\abs{t}^{\alpha} +o(|t|^{\alpha} ) 
\end{align*}
as $t\to 0$ for some $c$, $0<c<\infty$ and some $\alpha$, $0<\alpha\le 2$,   then, 
for any $\gamma>0$, $\frac{\xi_t}{t^{\gamma}}$ converges almost surely to zero as $t\to \infty$.  For the case of the fractional Brownian motion noise, the strong consistency is obtained from applying that theorem to the stationary Gaussian process
 \begin{align*}
 \xi_t:=\int_{-\infty}^{t} e^{-\theta (t-s)}\dif B^{H}_s
 \end{align*} in \cite{hu Nua 10, hu nua zhou 19}.
For our case, $\xi_t$ is not a stationary Gaussian process any more and thus we can not apply the above theorem of Pickands. We get around this difficulty by means of combining three techniques together: the hypercontractivity of multiple Wiener-It\^o integrals, Kolmogorov's continuity theorem and the Garsia-Rodemich-Rumsey inequality. Then we obtain the similar stochastic integrals (see Propositions~\ref{Ft 2 norm 1}  and \ref{prop ft ht}  below) converge almost surely to zero as $t$ tends to infinity. This method has been used to replace the  theorem of Pickands 
 in the literatures such as \cite{CZ 21, cai xiao 18, chw 17, shen 21, clp 21}. We would also like to mention \cite{msh} where they use the moment generation function method to obtain a similar result under some assumptions different to ours. Those two methods are similar since for multiple Wiener-It\^o integrals, the hypercontractivity can be used to obtain the exponential integrability, see e.g. \cite{Hu 07}. Please refer to \cite{Hu 07, St Va 79}  for the Garsia-Rodemich-Rumsey inequality.

\begin{thm}\label{main thm 2}
Assume $H\in (0,\,\frac38)$ and Hypothesis~\ref{hypthe 1} is satisfied. Then, both $\sqrt{T}( \hat{\theta}_T-\theta )$ and $\sqrt{T}( \tilde{\theta}_T-\theta ) $ are asymptotically normal as $T\to \infty$. Namely,
\begin{align}
\sqrt{T}( \hat{\theta}_T-\theta )& \stackrel{ {law}}{\to}  { \mathcal{N}(0, \,\theta \sigma_{H}^2)},\label{asy norm LSE}\\
\sqrt{T}( \tilde{\theta}_T-\theta )&\stackrel{ {law}}{\to}  \mathcal{N}(0, \,\theta \sigma_{H}^2 /4H^2),\label{asy norm SME}
\end{align}where
\begin{align*}
\sigma^2_H= (4H-1) +   \frac{2 \Gamma(2-4H)\Gamma(4H)}{\Gamma(2H)\Gamma(1-2H)}.
\end{align*}
\end{thm}
\begin{remark}
 Because of the calculation complication of the inequality \eqref{zhou ineq} in the case of $H\in (0,\frac12)$ (see Proposition~\ref{contraction ft}), we only obtain the CLT for $H\in (0,\,\frac38)$ in the present paper.  A more sharp estimate of the inequality \eqref{zhou ineq} is still awaiting. We will investigate this topic in other works. \end{remark}
\begin{thm}\label{B-E bound thm}
Let $\Phi(z)$ be the standard normal distribution function. 
Assume $H\in (0,\,\frac38)$ and Hypothesis~\ref{hypthe 1} is satisfied. Then, when  $H\neq \frac14$,  there exists a constant $C_{\theta, H} > 0$ such that when $T$ is large enough,
\begin{equation}\label{b-e bound 34}
\sup_{z\in \Rnum}\abs{P(\sqrt{\frac{T}{\theta \sigma^2_{H}}} (\hat{\theta}_T-\theta )\le z)-\Phi( z)}\le\frac{ C_{\theta, H}}{T^{ \delta}},
\end{equation}and
\begin{equation}\label{b-e bound 44}
\sup_{z\in \Rnum}\abs{P(\sqrt{\frac{4H^2 T}{\theta \sigma^2_{H}}} (\tilde{\theta}_T-\theta )\le z)-\Phi( z)}\le\frac{ C_{\theta, H}}{{T^{\delta} }},
\end{equation}
 where 
\begin{equation}\label{delta dingyi}
\delta=\left\{
      \begin{array}{ll}
 \frac12, & \quad \text{when } H\in {(0,\frac14)},\\ 
 {\frac32-4H}, &\quad \text{when } H\in (\frac14,\frac38).
 \end{array}
\right.
\end{equation}and when $H=\frac14$, the upper bound can be replaced by $\frac{ \log T}{\sqrt{T} }$.
\end{thm}
\begin{remark}
For the moment estimator $\tilde{\theta}_T$,  when the driven noise $G_t$ degenerates to the fractional Brownian $B_t^H$ with $H\in (0,\,\frac12)$, the Berry-Ess\'{e}en upper bound in \eqref{b-e bound 44} is shown to be also $\frac{1}{\sqrt{T}}$, see Proposition 4.1 of \cite{SV 18}. Thus, it is reasonable to conjecture that when $H\in (\frac14,\frac12)$, a better bound should also be $\frac{1}{\sqrt{T}}$ for the two estimators. This requires us to find two more sharp estimates than those given by Lemma 3.11 of \cite{CL 21} and Proposition~\ref{contraction ft} respectively.
\end{remark}

 In the remaining part of this paper, $C, c$ will be a generic positive constant independent of $T$ whose value may differ from line to line.

\section{Preliminary}\label{sec pre}
Denote $G = \set{G_t, t\in [0,T ]}$ as a continuous centered Gaussian process with covariance function
    $$ \E(G_tG_s)=R(s,t), \ s, t \in [0,T], $$
defined on a complete probability space $(\Omega, \mathcal{F}, P)$. The filtration $\mathcal{F}$ is generated by the Gaussian family $G$. Suppose in addition that the covariance function $R$ is continuous.
Let $\mathcal{E}$ denote the space of all real valued step functions on $[0,T]$. The Hilbert space $\mathfrak{H}$ is defined
as the closure of $\mathcal{E}$ endowed with the inner product
\begin{align*}
\innp{\mathbbm{1}_{[a,b)},\,\mathbbm{1}_{[c,d)}}_{\FH}=\E\big(( G_b-G_a) ( G_d-G_c) \big).
\end{align*}
Abusing the notation slightly, we also write $G=\{G(h), h \in \mathfrak{H}\}$ as the isonormal Gaussian process on the probability space $(\Omega, \mathcal{F}, P)$, indexed by the elements in the Hilbert space $\mathfrak{H}$. In other words, $G$ is a Gaussian family of random variables such that
   $$\mathbb{E}(G) = \mathbb{E}(G(h)) = 0, \quad \mathbb{E}(G(g)G(h)) = \langle g, h \rangle_{\mathfrak{H}} \,,$$
for any $g, h \in \mathfrak{H}$.

In particular, when $G $ is exact the fractional Brownian motion $B^H$, we denote by $\mathfrak{H}_1$ the associated Hilbert space.
\begin{nott}\label{innp h1 h2}
Denote by $\mathcal{V}_{[0,T]}$ the set of bounded variation functions on $[0,T]$. For  functions $f,\,g \in\mathcal{V}_{[0,T]}$, we define two  products as
\begin{align}
\innp{f,\,g }_{\mathfrak{H}_1}&=- \int_{[0,T]^{2}}f(t ) \frac{\partial R^B(t,s)}{\partial t}  \dif t \nu_{g}({\dif s}), \label{inner product 0010}\\
\innp{f,\,g }_{\mathfrak{H}_2}&=C_{H}^{'} \int_{[0,T]^{2}} \abs{f(t) g(s)}(ts)^{H-1}\dif t  \dif s =C_{H}^{'}  \mu(\abs{f})\mu(\abs{g}),\nonumber
\end{align} where ${\nu_{g}}$ is given below and $R^B(t,s)$ is the covariance function of the fractional Brownian motion as \eqref{exmp1} and $\mu(f)=\int f\dif\mu$ with $\mu(\dif x)=x^{\beta-1}\dif x$. We also denote for any $f\in L^2(R^2, \mu\times \mu)$
\begin{align*}
(\mu\times \mu) (f)= \int f \,\dif(\mu\times \mu).
\end{align*}
\end{nott}
The following proposition is an extension of \cite[Theorem 2.3]{Jolis2007} and \cite[Proposition 2.2]{CL 21}, which gives the inner products representation of the Hilbert space $\mathfrak{H}$:
\begin{proposition}\label{inner products represent}
$\mathcal{V}_{[0,T]}$ is dense in $\mathfrak{H}$ and we have
  \begin{align}  \label{inner product 001}
\innp{f,g}_{\FH}=\int_{[0,T]^2} R(t,s) \nu_f( \dif t)   \nu_{g}( \dif s),\qquad \forall f,\, g\in \mathcal{V}_{[0,T]},
\end{align}where $\nu_{g}$ is the restriction to $([0,T ], \mathcal{B}([0,T ]))$ of
the Lebesgue-Stieljes signed measure associated with $g^0$ defined as
\begin{equation*}
g^0(x)=\left\{
      \begin{array}{ll}
 g(x), & \quad \text{if } x\in [0,T),\\
0, &\quad \text{otherwise}.
 \end{array}
\right.
\end{equation*}
Furthermore, if the covariance function $R(t,s)$ satisfies Hypothesis~\ref{hypthe 1}, then
  \begin{align} \label{innp fg3}
\innp{f,g}_{\FH}=-\int_{[0,T]^2}  f(t) \frac{\partial R(t,s)}{\partial t} \dif t  \nu_{g}(\dif s),\qquad \forall f,\, g\in \mathcal{V}_{[0,T]},
\end{align}
and
\begin{equation}\label{inequality 29}
\abs{\innp{f,\,g}_{\mathfrak{H}} - \innp{f,\,g }_{\mathfrak{H}_1}}\leq  {\innp{f,\,g }_{\mathfrak{H}_2} },\qquad \forall f,\, g\in \mathcal{V}_{[0,T]}.
\end{equation}
\end{proposition}
\begin{corollary}\label{cor new}
 Let $g=h\cdot \mathbb{1}_{[a,b)}(\cdot)$ with $h$ a continuously differentiable function.  Then we have \cite{CL 21}
\begin{align*} 
\innp{g,\,g}_{\FH_1}&=-\int_{[a,b)^2}  h(t)h'(s) \frac{\partial R^B(t,s)}{\partial t} \dif t \dif s +\int_{[a,b)}  h(t)  \Big(h(b) \frac{\partial R^B(t,b)}{\partial t}-h(a) \frac{\partial R^B(t,a)}{\partial t} \Big)  \dif t.
\end{align*}
\end{corollary}
\begin{remark}\label{rem dandiao}
When $H\in (\frac12,\,1)$, the identity \eqref{inner product 0010}  is also equal to the following formula:
\begin{equation*}
\langle f,g\rangle_\mathfrak{H_1}=H(2H-1)\int_{[0,T]^{2}}f(t)g(s) \abs{t-s}^{2H-2}dtds,          \quad        \forall f,g\in \mathcal{V}_{[0,T]},
\end{equation*} and hence if $0\le f\le g$ then the monotonicity of the norm holds:
\begin{align*}
 \innp{f,\,f}_{\FH_1}\le \innp{g,\,g}_{\FH_1}  .
 \end{align*}
 A remarkable fact is that this monotonicity of the norm may not hold in the case of $H\in (0,\,\frac12)$. This is one of the reasons why it is more difficult to deal with the problems of $H\in (0,\,\frac12)$.  
\end{remark}
The following proposition is an immediately consequence of the identity (\ref{innp fg3}).
\begin{prop}
Suppose that Hypothesis~\ref{hypthe 1} holds. Then for any $ \varphi, \psi \in (\mathcal{V}_{[0,T]})^{\odot 2}$,
 \begin{align}
	\abs{\innp{\varphi,\, \psi}_{\FH^{\otimes 2}}  - \innp{\varphi,\,\psi}_{\FH_1^{\otimes 2}} } & \leq  (C_{H}' )^2(\mu \times \mu)(\abs{\varphi}) (\mu \times \mu)(\abs{\psi})  + 2C_{H}'(\mu \times \mu)(\abs{\varphi\otimes_{1'} \psi}), \label{inner product 00.ineq}
 \end{align} where $\varphi\otimes_{1'} \psi$ is the 1-th contraction between $\varphi$ and $\psi$ in $\mathfrak{H}_1^{\otimes 2}$, see \eqref{contractiondefn} below.
\end{prop}


Denote $\mathfrak{H}^{\otimes p}$ and $\mathfrak{H}^{\odot p}$ as the $p$th tensor product and the $p$th symmetric tensor product of the Hilbert space $\mathfrak{H}$. Let $\mathcal{H}_p$ be the $p$th Wiener chaos with respect to $G$. It is defined as the closed linear subspace of $L^2(\Omega)$ generated by the random variables $\{H_p(G(h)): h \in \mathfrak{H}, \ \|h\|_{\mathfrak{H}} = 1\}$, where $H_p$ is the $p$th Hermite polynomial defined by
$$H_p(x)=\frac{(-1)^p}{p!} e^{\frac{x^2}{2}} \frac{d^p}{dx^p} e^{-\frac{x^2}{2}}, \quad p \geq 1,$$
and $H_0(x)=1$. We have the identity $I_p(h^{\otimes p})=H_p(G(h))$ for any $h \in \mathfrak{H}$ where $I_p(\cdot)$ is the generalized Wiener-It\^o stochastic integral. Then the map $I_p$ provides a linear isometry between $\mathfrak{H}^{\odot p}$ (equipped with the norm $\frac{1}{\sqrt{p!}}\|\cdot\|_{\mathfrak{H}^{\otimes p}}$) and $\mathcal{H}_p$. Here $\mathcal{H}_0 = \mathbb{R}$ and $I_0(x)=x$ by convention.

We choose $\{e_k, k \geq 1\}$ to be a complete orthonormal system in the Hilbert space $\mathfrak{H}$. Given $f \in \mathfrak{H}^{\odot m}, g \in \mathfrak{H}^{\odot n}$, the $q$-th contraction between $f$ and $g$ is an element in $\mathfrak{H}^{\otimes (m+n-2q)}$ that is defined by
\begin{equation}\label{contractiondefn}
 f \otimes_q g = \sum_{i_1,\dots,i_q=1}^{\infty} \langle f,e_{i_1} \otimes \cdots \otimes e_{i_q} \rangle_{\mathfrak{H}^{\otimes q}} \otimes \langle g,e_{i_1} \otimes  \cdots  \otimes e_{i_q} \rangle_{\mathfrak{H}^{\otimes q}} \,,
\end{equation}
for $q=1,\dots, m \wedge n$. Especially, when $G $ is exact the fractional Brownian motion $B^H$, we denote by $f \otimes_q' g $ the $q$-th contraction between $f$ and $g$ is an element in $\mathfrak{H}_1^{\otimes (m+n-2q)}$ .
 
For $g \in \mathfrak{H}^{\odot p}$ and $h \in \mathfrak{H}^{\odot q}$, we have the following product formula for the multiple integrals,
\begin{equation}\label{ito.prod}
	I_p(g) I_q(h) = \sum_{r=0}^{p \wedge q} r!\binom{p}{r} \binom{q}{r} I_{p+q-2r}(g\tilde\otimes_r h) \,,
\end{equation}
where $g\tilde\otimes_r h$ is the symmetrization of $g\otimes_r h$ (see \cite{Nou 12}).

The following Theorem \ref{fm.theorem}, known as the fourth moment theorem, provides necessary and sufficient conditions for the  convergence of a sequence of random variables to a normal distribution (see \cite{Nou 12, Nua Pec 05}).
\begin{thm} \label{fm.theorem}
 Let $ p \geq 2 $ be a fixed integer. Consider  a collection of elements $\{f_{T}, T>0\}$ such that $f_{T} \in \mathfrak{H}^{\odot p}$ for every $T>0$. Assume further that
 \[
  \lim_{T \to \infty}\mathbb{E} [I_p (f_T) ^2] = \lim_{T \to \infty} p! \|f_T\|^2_{\mathfrak{H}^{\otimes{p}}} = \sigma^2 .
 \]
Then the following conditions are equivalent:
 \begin{enumerate}
	\item $\lim_{T \to \infty} \mathbb{E}[I_p(f_T)^4] = { 3\sigma^4}$.
	\item For every $q=1,\dots, p-1$, $\lim_{T \to \infty}||f_T \otimes_q f_{T}||_{\mathfrak{H}^{\otimes 2(p-q)}} = 0$.
    \item As $T $ tends to infinity, the $p$-th multiple integrals $\{I_p (f_T), T \geq 0\}$ converge in distribution to a Gaussian random variable $N(0,\sigma^2)$.
 \end{enumerate}	
\end{thm}

\medskip
\section{Strong Consistency:  Proof of Theorem~\ref{thm strong} }
 We will discuss exclusively the case $H \in(0,\, \frac12)$ since the case of $H> \frac12$ has been treated in \cite{CZ 21}.
We first define some important functions that will be used in the proof. Denote
\begin{align}
f_T(t,s)&= e^{-\theta \abs{t-s}}\mathbbm{1}_{\set{0\le s,t\le T}},\label{ft ts 000} \\
h_T(t,s)&= e^{-\theta (T-t)-\theta (T-s)}\mathbbm{1}_{\set{0\le s,t\le T}},\label{ht ts}\\
g_T(t,s)&= \frac{1}{2\theta T}(f_T- h_T).\label{gt ts}
\end{align}
The solution to the SDE \eqref{fOU} with $X_0=0, \sigma=1$ is $$X_t = \int_0^t e^{-\theta(t-s)} dG_s = I_1(f_T(t, \cdot) \mathbbm{1}_{[0,t]}(\cdot))\,.$$
We apply the product formula of multiple integrals \eqref{ito.prod} and stochastic Fubini theorem to obtain
\begin{align}
\frac{1}{T} \int_0^T X_t^2\mathrm{d} t =  I_2(g_T)+b_T,  \label{x square}
\end{align}
where
 \begin{equation}
	 b_T=\frac{1}{T}\int_0^T\, \norm{e^{-\theta (t-\cdot) }\mathbbm{1}_{[0,t]}(\cdot)}^2 _{\mathfrak{H}}\dif t.\label{bt bt}
 \end{equation}
From the equation (\ref{hattheta}), we can write
\begin{align} \label{ratio 1}
\sqrt{T} (\hat{\theta}_T-\theta )&=-\frac{\frac{1}{2\sqrt{T}} I_2(f_T)}{ I_2(g_T)+b_T}.
\end{align}

For the kernel $f_T(r,s)\mathbbm{1}_{\set{0\le r,s\le t}}$, it's double Wiener-It\^{o} integral  
\begin{equation}\label{Ft}
	F_t:=I_2(f_T(r,s)\mathbbm{1}_{\set{0\le r,s\le t}}),\quad t\in [0,T]
\end{equation} is a stochastic process which is named as a chaos (stochastic) process \cite{HPAU 94}.
The next two propositions are about the asymptotic behaviors of the second moment of $F_t$ and the increment $F_t - F_s$ respectively, which are used to estimate the modulus of continuity of the chaos process $\set{F_t}$ on $[n,\,n+1]$ for any integer $n\ge 1$ and then to get the asymptotic growth of the process $F_t$ as $t\to \infty$ (see Proposition~\ref{prop ft ht}). Those types of sample path properties of chaos processes are the key tools to show the strong consistency in the present paper.
\begin{prop}\label{Ft 2 norm}
When $H\in (0,\,\frac12)$,
\begin{align}\label{Ft-norm}
\lim_{t\to\infty}\frac{1}{4 \theta \sigma_{H}^2 t}\E[\abs{F_t} ^2]= (H   \Gamma(2H ) \theta^{-2H})^2 .
\end{align}
\end{prop}
\begin{proof}
The inequality (\ref{inner product 00.ineq}) implies that
\begin{align}\label{ft norm 02}
\abs{ \norm{f_t}^2_{\FH^{\otimes 2}}-\norm{ f_t}_{\FH_1^{\otimes 2}}^2} \le  \big(C_{H}'(\mu\times \mu)(f_t) \big)^2 +2C_{H}'(\mu\times \mu)(f_t\otimes_{1'}f_t).
\end{align}
 Lemma 17 in \cite{hu nua zhou 19}  implies that when $H\in (0,\,\frac12)$,
\begin{align}\label{fth1 zhouh}
\lim_{T\to \infty } \frac{1}{2 \theta \sigma_{H}^2 t}\norm{ f_t}_{\FH_1^{\otimes 2}}^2&=(H   \Gamma(2H )\theta^{-2H}  )^2.
\end{align}
 Lemma~\ref{upper bound F} implies that there exists a positive constant $C$ independent on $T$ such that 
\begin{align}\label{ft fh2 norm 2}
 (\mu\times \mu)(f_t) =\int_{[0,t]^2} e^{-\theta\abs{r-s}}(rs)^{H-1}\dif r\dif s=2\int_0^t  e^{-\theta r }r^{H-1}\dif r\int_0^r e^{ \theta s }s^{H-1}\dif s  \le C .
\end{align}
Then it follows from Lemma~\ref{lem 3.2 guji ft yasuo} that when $H\in (0,\,\frac12)$,
	$$\lim_{t\to\infty}\frac{1}{  t} \abs{\norm{f_t}^2_{\FH^{\otimes 2}}-\norm{f_t}_{\FH_1^{\otimes 2}}^2} =0 , $$
which together with \eqref{fth1 zhouh} and the {It\^{o}'s isometry} $$\E[\abs{F_t}^2]=2\norm{f_t}^2_{\FH^{\otimes 2}}$$ implies the desired (\ref{Ft-norm}).
\end{proof}
\begin{remark}\label{remark ft2 shoulian rate}
Together with Lemma 3.11 of \cite{CL 21}, we have that when $H\in (0,\frac12)$, 
the speed of convergence  $$\frac{1}{2 \theta \sigma_{H}^2 t}\norm{f_t}^2_{\FH^{\otimes 2}} \to (H   \Gamma(2H ) \theta^{-2H})^2 $$ is  at least $\frac{1}{t^{1-2H}}$ as $t\to \infty$.
\end{remark}
\begin{nott}Let $0\le s<t\le T$. Denote 
 \begin{align}
\phi_1(u,v)&=e^{-\theta \abs{u-v}}\mathbbm{1}_{\set{s\le u,v\le t}},\label{pjhi1uv}\\
\phi_2(u,v)&=e^{-\theta \abs{u-v}}(\mathbbm{1}_{\set{0\le u\le s,\, s\le v\le t}} + \mathbbm{1}_{\set{0\le v\le s,\, s\le u\le t}} ).\label{pjhi12 uv}
\end{align}
\end{nott}
\begin{prop}\label{Ft 2 norm 1} 
If Hypothesis~\ref{hypthe 1} is satisfied,
there exists a  constant $C>0$ independent of $T$ such that for all $s,t\ge 0$, 
\begin{align}\label{E norm 2}
\E[\abs{F_t-F_s} ^2]\le C \Big( \abs{t-s}^{4H+2}+ \abs{t-s}^{4H+1}+\abs{t-s}^{4H} +\abs{t-s}^{2H+1} +\abs{t-s}^{2H}+\abs{t-s}^{H}\Big).
\end{align}
Moreover, for any real number $p > \frac{4}{H}, q >1$ and integer $n \ge1$,  there exists a random constant $R_{p,q}$ independent of $n$  such that
\begin{align}\label{zhengze cha ft}
\abs{F_t-F_s}\le R_{p,q} n^{q/p} , \qquad \forall \ t,s\in [n,n+1].
\end{align}
\end{prop}
\begin{proof} 
 It\^{o}'s isometry implies that
\begin{align}\label{ft fs}
\E[\abs{F_t-F_s} ^2]&=2\norm{f_t-f_s}^2_{\FH^{\otimes 2}} \le 4 (\norm{\phi_1}_{\FH^{\otimes 2}}^2+\norm{\phi_2 }_{\FH^{\otimes 2}}^2),
\end{align}
 For simplicity, we can assume that $\theta=1$.
 Lemma~\ref{upper bound F} implies that there exists a positive constant $C$ independent on $T$ such that 
\begin{align*}
(\mu\times \mu)( \abs{\phi_1}) &=2\int_s^t e^{-u}u^{H-1}\dif u\int_s^u e^{v}v^{H-1}\dif v\le 2\int_s^t  u^{H-1}\dif u\le C \abs{t-s}^{H},\\
(\mu\times \mu)( \abs{\phi_2})&\le 
 2\int_0^s e^{-  (s-u)}u^{H-1}\dif u \int_s^t e^{-  (v-s)}v^{H-1}\dif v\le 2\int_s^t  v^{H-1}\dif v\le  C\abs{t-s}^H.
\end{align*}
Hence,  it follows from the inequality \eqref{inner product 00.ineq} and Lemma~\ref{lem phi guji}  that
\begin{align}\label{phi1 norm}
\norm{ \phi_1}_{\FH^{\otimes 2}}^2&\le \norm{ \phi_1}_{\FH_1^{\otimes 2}}^2 + \big(C_H' (\mu\times \mu)( \abs{\phi_1})\big)^2
+2C_H'(\mu\times\mu)(\phi_1\otimes_{1'}\phi_1)\nonumber \\
&\le C\Big(\abs{t-s}^{2H}+\abs{t-s}^{4H}+\abs{t-s}^{4H+1}+\abs{t-s}^{4H+2}\Big),
\end{align}
 and similarly, the inequality \eqref{inner product 00.ineq} and Lemma~\ref{lem phi2 guji}  implies that
\begin{align}\label{phi2 norm}
\norm{ \phi_2}_{\FH^{\otimes 2}}^2& \le \norm{ \phi_2}_{\FH_1^{\otimes 2}}^2 + \big(C_H' (\mu\times \mu)( \abs{\phi_2})\big)^2
+2C_H'(\mu\times\mu)(\phi_2\otimes_{1'}\phi_2)\nonumber \\
&\le C \big(\abs{t-s}^{2H+1} +\abs{t-s}^{2H}+\abs{t-s}^{H}\big).
\end{align}
By plugging the inequalities \eqref{phi1 norm} and \eqref{phi2 norm} into \eqref{ft fs}, we obtain the desired estimate \eqref{E norm 2}.

The inequality \eqref{zhengze cha ft} is then obtained from the Garsia-Rodemich-Rumsey inequality (see Proposition 3.4 of \cite{chw 17}). For the reader's  convenience, we rewrite it as follows:  By the inequality \eqref{E norm 2}, there exists a positive constant $C$ independent of $T$ such that, for any $\abs{t-s}\le 1$, 
\begin{align*}
\E[\abs{F_t-F_s} ^2]\le C \abs{t-s}^{H}.
\end{align*}
The hypercontractivity of multiple Wiener-It\^o integrals implies that for any $p\ge 2$ and any $n\le t<s\le n+1$,
\begin{align*}
\E[\abs{F_t-F_s} ^p]\le C\times  (p-1)^p\abs{t-s}^{pH/2}.
\end{align*} 
Then $F_t$ has a continuous realization on $[n,n+1]$ for all integer $n\ge 1$ by Kolmogorov's continuity theorem. 
Next, take $\Psi(x)=x^p$ and $\rho(x)=x^{H/2}$. Denote 
$$B_n=\int_{[n,n+1]^2 }\Psi\big(\frac{\abs{F_t-F_s}}{\rho(\abs{t-s})}\big) \dif t \dif s.$$
The above inequality implies that for any $q>1$,
\begin{align*}
\E\big( \sum_{n=1}^{\infty}  \frac{B_n}{n^q} \big)=\sum_{n=1}^{\infty}  \frac{\E (B_n)}{n^q}\le C\times  (p-1)^p \sum_{n=1}^{\infty}  \frac{1}{n^q}<\infty.
\end{align*}
Hence, there exists a random constant $R_{p,q}$ such that 
$$\sum_{n=1}^{\infty}  \frac{B_n}{n^q} \le R_{p,q},$$
which implies that for all positive $q>1$ and integer $n\ge 1$,
\begin{equation}\label{bn bound uniform}
 {B_n}\le  R_{p,q}{n^q}.\end{equation}
An application of the Garsia-Rodemich-Rumsey inequality,  see e.g. \cite{ Hu 07, St Va 79}, implies that  when $p>\frac{4}{H}$,  
 we have  for all $s,t\in  [n,n+1]$,
 \begin{align*}
\abs{F_t-F_s}\le  8\int_0^{\abs{t-s}}\Psi^{-1}(\frac{4B_n}{u^2})\rho'(u)\dif u= 2H\frac{(4B_n)^{\frac{1}{p}}}{H-\frac{4}{p}}\abs{t-s}^{\frac{H}{2}-\frac{2}{p}}<c_p B_n^{\frac{1}{p}}.
\end{align*} This combined with \eqref{bn bound uniform} 
shows the proposition.
\end{proof}
We denote a chaos process 
\begin{equation}\label{ht}
	J_t=I_2(h_T(r,u)\mathbbm{1}_{\set{0\le r,u\le t}}),\quad t\in [0,T],
\end{equation}
and apply the similar computations as above to obtain its  sample path properties as follows: 

\begin{prop}\label{prop ht} If Hypothesis~\ref{hypthe 1} is satisfied, there exists a constant $C>0$ independent of $T$ such that
\begin{align}\label{ht Jt upper bounds}
\sup_{t\ge 0}\E[\abs{J_t} ^2]&< C,
\end{align}
and there exist two positive constants $C$ and $\alpha\in (0,1)$ independent of $T$ such that, for any $\abs{t-s}\le 1$,
\begin{align*}
\E[\abs{J_t-J_s} ^2]&\le C \abs{t-s}^{\alpha}.
\end{align*}
Moreover, for any real number {$p > \frac{2}{\alpha}, q >1$} and integer $n \ge1$,  there exists a random constant $R_{p,q}$ independent of $n$  such that
\begin{align}
\abs{J_t-J_s}\le R_{p,q} n^{q/p} , \qquad \forall \ t,s\in [n,n+1].
\end{align}
\end{prop}
\begin{prop}\label{prop ft ht} 
Let $F_T$ and $J_T$ be given in \eqref{Ft} and \eqref{ht} respectively. If Hypothesis~\ref{hypthe 1} is satisfied,
we have
\begin{align}
\lim_{T\to \infty}\frac{F_T}{T}=0 \quad {\rm and} \quad \lim_{T\to \infty}\frac{J_T}{ {T^{\alpha}}}=0, \qquad a.s.
\end{align} for any $\alpha>0$.
\end{prop}
\begin{proof}
The proof is similar to that in \cite{CZ 21} or \cite{chw 17}. We will only show $\lim\limits_{T\to \infty}\frac{F_T}{T}=0$, and the other is similar. When $H \in (0,\frac{1}{2})$, Chebyshev's inequality, the hypercontractivity of multiple Wiener-It\^{o} integrals and Proposition~\ref{Ft 2 norm} imply that for any $\epsilon > 0$,
\begin{equation*}P\left(\frac{\abs{F_n}}{n} > \epsilon \right) \leq \frac{\mathbb{E}F_n^4}{n^4 \epsilon^4} \leq \frac{C\left(\mathbb{E}(F_n^2)\right)^2}{n^4 \epsilon^4}
 \leq C n^{-2}  \,,
 \end{equation*}which implies that $\frac{F_n}{n}$ converges to $0$ almost surely as $n\to \infty$ by the Borel-Cantelli lemma.\\
Since
\begin{align*}
\abs{\frac{F_T}{T}}\le \frac{1}{T}\abs{F_T-F_n}+ \frac{n}{T}\frac{\abs{F_n}}{n},
\end{align*}where $n=[T]$ is the biggest integer less than or equal to a real number $T$, we have $\frac{F_T}{T}$ converges to $0$ almost surely as $T\to \infty$ by Proposition~\ref{Ft 2 norm 1}.
\end{proof}

Proposition \ref{prop ft ht} implies the following chaos process 
$$I_2(g_T) = \frac{F_T - J_T}{2\theta T} \to 0$$ as $T \to \infty$ almost surely. Next we study the term $b_T$.
\begin{prop}\label{prop bt lim}
Let $b_T$ be given by \eqref{bt bt}. Suppose that Hypothesis~\ref{hypthe 1} holds. We have
\begin{align}\label{prop bt limit}
\lim_{T\to\infty}b_T = H\Gamma(2H)\theta^{-2H} >0.
\end{align}
\end{prop}
\begin{proof}
The L'H\^{o}pital's rule implies that
\begin{align}\label{bt limit}
\lim_{T\to\infty}b_T&= \lim_{T\to\infty} \norm{e^{-\theta (T-\cdot) }\mathbbm{1}_{[0,T]}(\cdot)}^2 _{\mathfrak{H}}.
\end{align}It follows from  from \cite{hu nua zhou 19} or \cite{CL 21} that

\begin{align*}
\lim_{T\to\infty} \norm{e^{-\theta (T-\cdot) }\mathbbm{1}_{[0,T]}(\cdot)}^2 _{\mathfrak{H}_1}=H\Gamma(2H)\theta^{-2H}.
\end{align*}
The identity (\ref{inequality 29}) implies that when $T\ge 1$,
\begin{align}
\abs{ \norm{e^{-\theta (T-\cdot) }\mathbbm{1}_{[0,T]}(\cdot)}^2 _{\mathfrak{H}}- \norm{e^{-\theta (T-\cdot) }\mathbbm{1}_{[0,T]}(\cdot)}^2 _{\mathfrak{H}_1}}&\le C\abs{\int_0^T e^{-\theta (T-u)} u^{H-1} \dif u}^2 \le C T^{2(H-1)}. \label{sudu}
\end{align}
 where in the last line we use Lemma~\ref{upper bound F}. By substituting the above limit and estimate into the identity \eqref{bt limit}, we obtain the desired limit \eqref{prop bt limit}.
\end{proof}
\begin{remark} \label{rem 36}
 Lemma 3.2 of \cite{CL 21} implies that the speed of convergence
$$ \frac{1}{T}\int_0^T  \norm{e^{-\theta (t-\cdot) }\mathbbm{1}_{[0,t]}(\cdot)}^2 _{\mathfrak{H}_1}  \dif t \to  H\Gamma(2H)\theta^{-2H} $$ is at least $\frac{1}{T}$. Clearly, $2(H-1)<-1$ when $H\in (0,\frac12)$. Then by the inequality (\ref{sudu}), there exists a positive constant $C$ such that when $T$  is large enough,
\begin{equation}\label{bt speed}
\abs{b_T- H\Gamma(2H)\theta^{-2H}} \leq \frac{C}{T } .\end{equation}
\end{remark}

\noindent{\it Proof of Theorem~\ref{thm strong}.\,} From\eqref{x square},  \eqref{gt ts},  \eqref{Ft}, and \eqref{ht},
\begin{align*}
\frac{1}{T} \int_0^T X_t^2\mathrm{d} t&=  I_2(g_T)+b_T=\frac{1}{2\theta}[\frac{F_T}{T} - \frac{J_T}{T}]+ b_T .
\end{align*}
Proposition~\ref{prop ft ht} and \ref{prop bt lim} imply that
\begin{equation*}
\lim_{T\to \infty}\frac{1}{T} \int_0^T X_t^2\mathrm{d} t =  H\Gamma(2H)\theta^{-2H},\,\,\, a.s.,
\end{equation*}which implies that the moment estimator $\tilde{\theta}_T$ is strongly consistent by the continuous mapping theorem. Since \begin{align*} 
 \hat{\theta}_T-\theta &=\frac{- \frac{1}{2T} F_T}{ I_2(g_T)+b_T},
\end{align*}
 Proposition~\ref{prop ft ht} and \ref{prop bt lim} imply that the least squares estimator $\hat{\theta}_T$ is also strongly consistent.
{\hfill\large{$\Box$}}
\section{{The Asymptotic normality and the Berry-Ess\'{e}en bound}}
\subsection{The Asymptotic normality}

\begin{prop}\label{contraction ft} Let $\delta $ be given as in  \eqref{delta dingyi}. 
When $H\in (0,\,\frac14)\cup(\frac14,\frac38)$, there exists a constant $C_{\theta,H} > 0$ such that
\begin{equation}\label{zhou ineq}
\frac{1}{T}\norm{f_T\otimes_{1} f_T}_{\mathfrak{H}^{\otimes 2}}\le
 \frac{C_{\theta,H}}{T^{\delta }},
\end{equation} and when $H=\frac14$, the upper bound can be replaced by $\frac{ \log T}{\sqrt{T} }$.
\end{prop}
\begin{proof} Without loss of generality, we assume $\theta=1$.
	Recall that
	\begin{align*}
	\left(f_T\otimes_1 f_T\right) (u_1,u_2) &:= -\int_{[0,T]^2} f_T(u_1,v_1) \frac{\partial}{\partial v_2}f_T(u_2,v_2)  \frac{\partial   }{\partial v_1}R(v_1,\,v_2) \dif v_1\dif v_2\\
	&=\int_{[0,T]^2}e^{-\abs{u_1-v_1}-\abs{u_2-v_2}}\mathrm{sgn}(v_2-u_2) \frac{\partial   }{\partial v_1}R(v_1,\,v_2)  \dif v_1\dif v_2\\
	&+ \int_{[0,T]}e^{-\abs{u_1-v_1}- (T-u_2)}  \frac{\partial   }{\partial v_1}R(v_1,\,T)    \dif v_1.
	\end{align*}
The triangle inequality implies
\begin{equation}
\|f_T \otimes_1 f_T \|^2_{\mathfrak{H}^{\otimes 2}} \le    \norm{f_T \otimes_{1'} f_T}^2_{\mathfrak{H}^{\otimes 2}} +  \norm{f_T \otimes_1 f_T-f_T \otimes_{1'} f_T}^2_{\mathfrak{H}^{\otimes 2}}  .\label{ft contraction sum}
\end{equation}
It follows from the inequality (\ref{inner product 00.ineq}) that when $H\neq \frac14$,
\begin{align*}
 \norm{f_T \otimes_{1'} f_T}^2_{\mathfrak{H}^{\otimes 2}}& \le  \norm{f_T \otimes_{1'} f_T}^2_{\mathfrak{H}_1^{\otimes 2}}\\
 &+C \Big[\big((\mu\times \mu)(\abs{f_T \otimes_{1'} f_T})\big)^2+ (\mu\times \mu)\big(\abs{(f_T \otimes_{1'} f_T)\otimes_{1'} (f_T \otimes_{1'} f_T)}\big)\Big]\\
 &\le C [ T+T^{2\gamma_1} +T^{\gamma}],
\end{align*} where in the last line, we use the inequality (3.17) of  \cite{hu nua zhou 19} and Lemma~\ref{lem 3.2 guji ft yasuo} and Lemma~\ref{lemma last 2-kappa} respectively. 
The formula of integration by parts implies that
\begin{align*}
f_T \otimes_1 f_T-f_T \otimes_{1'} f_T= \int_{[0,T]^2} f_T(u_1,v_1) f_T(u_2,v_2) \frac{\partial ^2  }{\partial v_1 \partial v_2}\big(R(v_1,\,v_2) -R^B(v_1,\,v_2)\big)  \dif v_1\dif v_2,
\end{align*} which, together with Lemma~\ref{chennew 000}, implies that
\begin{equation}\label{zhou ineq chennew 000}
 \norm{f_T \otimes_1 f_T-f_T \otimes_{1'} f_T}^2_{\mathfrak{H}^{\otimes 2}}\le  {C_{\theta,H}}{T^{2\gamma_1}}.
\end{equation}

Comparing three values $1,\,2\gamma_1$, and $\gamma$, we see that the 
largest one is $$\delta_0=\max\set{1,\,2\gamma_1,\,\gamma}=\left\{
      \begin{array}{ll}
 1, & \quad \text{when } H\in (0,\frac14],\\ 
 {8H-1}, &\quad \text{when } H\in (\frac14,\frac12).
 \end{array}
\right.$$ Then we have 
\begin{align*} \norm{f_T \otimes_1 f_T }^2_{\mathfrak{H}^{\otimes 2}} \le  CT^{\delta_0}. \end{align*}
Clearly, $\delta=1-\frac{\delta_0}{2} $. Hence, we have the desired \eqref{zhou ineq}. The case of $H=\frac14$ is similar.
\end{proof}
\noindent{\it Proof of Theorem~\ref{main thm 2}.\,} Denote a constant that depends on $\theta$ and $H$ as $$a:=   H\Gamma(2H)\theta^{-2H}.$$
First, Proposition ~\ref{Ft 2 norm}, Proposition \ref{contraction ft} and Theorem \ref{fm.theorem} imply that when $H\in (0,\frac38)$, as $T\to\infty$,
\begin{equation} \label{asy norm Ft}
\frac{1}{2\sqrt{T}} F_T  \stackrel{ {law}}{\to} \mathcal{N}(0,\,\theta a^2\sigma_{H}^2).
\end{equation}
Second, recall the identity (\ref{ratio 1}):
\begin{align*}
\sqrt{T}( \hat{\theta}_T-\theta) =\frac{-\frac{1}{2\sqrt{T}} F_T}{ I_2(g_T)+b_T}.
\end{align*}
The Slutsky's theorem, Proposition~\ref{prop bt lim}  and the convergence result (\ref{asy norm Ft})  imply that the asymptotic normality (\ref{asy norm LSE}) holds when $H\in (0,\frac38)$.

Third,  we can show the asymptotic normality of
 \begin{align}\label{asy norm moment}
 \sqrt{T}\Big(\frac{1}{T} \int_0^T X_s^2\mathrm{d} s - a \Big)  \stackrel{ {law}}{\to} \mathcal{N}(0,\,a^2\sigma_{H}^2/\theta)
 \end{align}when $H\in (0,\frac38)$. 
 In fact, we have
 \begin{align}\label{asy norm moment 2}
 \sqrt{T}\Big(\frac{1}{T} \int_0^T X_s^2\mathrm{d} s -  a \Big)&= \frac{1}{2\theta}\Big[\frac{F_T}{\sqrt{T}} - \frac{J_T}{\sqrt{T}}\Big]+ \sqrt{T} \big(b_T-a) .
 \end{align}
The inequality \eqref{bt speed} implies that
$$
 \lim_{T\to\infty} \sqrt{T} \big(b_T-a)=0.$$
 Proposition~\ref{prop ht} and Proposition~\ref{prop ft ht} imply that $\frac{J_T}{\sqrt{T}}\to 0 $ a.s. as $T\to\infty$. Thus, the Slutsky's theorem implies that (\ref{asy norm moment}) holds. Finally, since
\begin{align*}
 \tilde{\theta}_T=\Big( \frac{1}{ {H}  \Gamma(2H )T} \int_0^T X_s^2\mathrm{d} s\Big)^{-\frac{1}{2H}} ,
\end{align*}
 the delta method  
 implies that the asymptotic normality (\ref{asy norm SME}) holds when $H\in (0,\frac38)$.
{\hfill\large{$\Box$}}\\

 \subsection{the Berry-Ess\'{e}en bound}
\begin{prop}\label{prop 4.2 be bound}
Let the constant $\delta $ be given as in \eqref{delta dingyi} 
and the double Wiener-It\^o integrals $F_T,\,J_T $ be as in \eqref{Ft} and \eqref{ht} respectively. Denote the double Wiener-It\^o integral
$$Q_T := \frac{F_T-J_T }{ \sqrt{T}}  .$$
Then there exists a constant $C_{\theta, H}$ such that when $H\neq \frac14$ and $T$ is large enough
 \begin{align}
 \sup_{z\in \Rnum}\abs{ P(\frac{1}{\sqrt{ 4\theta a^2\sigma_{H}^2}} Q_T\le z)- \Phi( z) }\le \frac{ C_{\theta, H}}{{T^{\delta}}},
 \end{align}where $\Phi(u)$ stands for the standard normal distribution function. When $H= \frac14$, the upper bound can be replaced by $\frac{\log T}{\sqrt{T}  }$.
\end{prop}
\begin{proof}We follow the line of the proof of Proposition 5.3  of \cite{CZ 21}. It follows from Remark~\ref{remark ft2 shoulian rate}, the inequality \eqref{ht Jt upper bounds} and Lemma~\ref{ft ht neiji bounds} that 
\begin{align*}
\frac{\abs{\mathbb{E}[Q_T^2]- {4}\theta a^2\sigma_{H}^2 }}{\mathbb{E}[Q_T^2] \vee ( {4}\theta a^2\sigma_{H}^2)}\le C\times (T^{2H-1}+T^{-1}+T^{(4H-1)_{+}-1}\log T)  \le C T^{2H-1}.
\end{align*} Clearly, $1-2H\ge \delta$. Then the the Fourth moment Berry-Ess\'{e}en bound (see, for example, Corollary 5.2.10 of \cite{Nou 12}) implies that we only need to show that when $T$ is large enough,
 \begin{align}\label{ Fourth moment Berry-Esseen bound}
 \E[Q_T^4]-3 \E[Q_T^2]^2\le \frac{C}{T^{ 2\delta  }},
 \end{align} which is the consequence of the following three estimates from the identity (5.8) of \cite{CZ 21}.\\
Proposition \ref{contraction ft} implies that 
  \begin{align*}
 \E\left[\left(\frac{F_T }{  \sqrt{T}}\right)^4\right] -3 \left[\E \left(\frac{F_T }{   \sqrt{T}}\right)^2\right]^2 \le C \left(\frac{1}{T}\|f_T \otimes_1 f_T\|_{\mathfrak{H}^{\otimes 2}}\right)^2\le \frac{C}{T^{2\delta}} \,.
 \end{align*}
The Cauchy-Schwarz inequality and Lemma~\ref{ft ht neiji bounds} imply that when $H\neq \frac14$ 
 \begin{align*}
\abs{ \left[\E[Q_T^2]\right]^2 - \left[\E \left(\frac{F_T }{   \sqrt{T}}\right)^2\right]^2 }&\le  \E\left[ Q_T^2 +\left(\frac{F_T }{   \sqrt{T}}\right)^2\right]   \E \abs{  \frac{H_T }{   \sqrt{T}} \left(\frac{H_T-2 F_T }{\sqrt{T}} \right)   }  \\
& \le {C}\times \frac{1}{ {T}^{1-(4H-1)_{+} }}  \le \frac{C}{T^{2\delta}}.
 \end{align*} 
 The identity (5.9) of \cite{CZ 21}, Lemma~\ref{ft ht neiji bounds} and Proposition \ref{contraction ft} imply that when $H\neq \frac14$ 
  \begin{align*}
\frac{1}{T^2}\abs{  \E[F_T^3 J_T] }&\le \frac{2}{T^2}\Big[  \norm{f_T {\otimes}_1f_T }_{\FH^{\otimes 2}} \norm{f_T}_{\FH^{\otimes 2}}  \norm{h_T}_{\FH^{\otimes 2}}  +\norm{f_T}_{\FH^{\otimes 2}}^2\abs{\innp{  f_T,\,  h_T}_{\mathfrak{H}^{\otimes 2}} }\Big]  \\
& \le \frac{C}{T^2}\times\big[ T^{\frac32-\delta}+T^{1+ (4H-1)_{+}} \big]  \le \frac{C}{T^{2\delta}}.
 \end{align*} 
\end{proof}
\noindent{\it Proof of Theorem~\ref{B-E bound thm}.\,} 
Denote $a=  H\Gamma(2H)\theta^{-2H}  $ and $b_T$ be given by \eqref{bt bt}. Then Remark~\ref{rem 36} and Remark~\ref{remark ft2 shoulian rate} imply that  there exists a constant $C>0$ such that  for $T$ large enough,
\begin{align*}
\abs{b_T^2-\frac{\norm{f_T}^2_{\mathfrak{H}^{\otimes 2}}}{2\theta \sigma^2_{H}T}}\le \abs{b_T^2- a^2}+\abs{\frac{\norm{f_T}^2}{2\theta \sigma^2_{\beta}T}-a^2}\le  C\times [\frac{1}{T}+ \frac{1}{T^{ 1-2H}}]\le \frac{C}{T^{1-2H}}.
\end{align*}
Recalling the definition of $g_T$ in \eqref{gt ts}, from Propositions~\ref{Ft 2 norm} and \ref{prop ht}, we have that there exists a constant $C>0$ such that  for $T$ large enough,
\begin{align*}
\norm{g_T}_{\FH^{\otimes 2}}\le \frac{C}{\sqrt{T}}.
\end{align*}
Hence, Corollary 1 of \cite{kim 3} implies that when $H\neq \frac14$, there exists a positive constant $C $ independent on $T$ such that, for $T$ large enough,
\begin{align*}
    &\sup_{z\in \Rnum}\abs{P\Big(\sqrt{\frac{T}{\theta \sigma^2_{H}} }(\hat{\theta}_T-\theta )\le z\Big)-\Phi(z)}   \\
& \le C \times \max \Big(\abs{b_T^2-\frac{\norm{f_T}_{\mathfrak{H}^{\otimes 2}}^2}{2\theta \sigma^2_{H}T}},\,\frac{1}{T}\norm{f_T\otimes_1 f_T} ,\, \norm{\frac{f_T}{\sqrt{T}}} \cdot\norm{g_T} ,\,  \norm{g_T}_{\FH^{\otimes 2}}^2 \Big) \\
&\le C\times [\frac{1}{T^{1-2H}}+\frac{1}{T^{\delta}}+ \frac{1}{\sqrt{T}}],
\end{align*}where in the last line we use Propositions~\ref{Ft 2 norm} and \ref{contraction ft}. Comparing three values $1/2,\, 1-2H$, and $\delta$, we see that the  smallest one is $\delta$. Hence, we obtain the Berry-Ess\'{e}en type bound \eqref{b-e bound 34} of the LSE $\hat{\theta}_T$.

The Berry-Ess\'{e}en type bound \eqref{b-e bound 44} of the ME $\,\tilde{\theta}_T$ can be obtained along the similar arguments of Theorem 3.2 in \cite{SV 18}. Please refer to the proof of Theorem 1.3 in \cite{CZ 21} for details. We only point out that \eqref{b-e bound 44} is a consequence of Proposition~\ref{prop 4.2 be bound} together with the following basic estimate:
 \begin{align}\label{zuihou gujishizi}
\sup_{z> -c \sqrt{T}}\Big[ \abs{ {\Phi} \left(\nu(z)- \frac{ { \sqrt{\theta T}} (b_T-a)}{\sqrt{   a^2\sigma_{H}^2}} \right)- {\Phi}(\nu(z)) } + \abs{\bar{\Phi}(\nu(z))-\Phi( z) }\Big]\le \frac{C}{\sqrt{T}},
 \end{align}
 where $c=\frac{2H\sqrt{\theta}}{\sigma_H}$  and $ \bar{\Phi}(z)=1-\Phi( z)$ and
 $$ \nu(z)=  \frac{c}{ 2H}  \sqrt{T}  \Big[ \big( 1+\frac{z }{c \sqrt{T}}\big)^{-2H}-1 \Big]. $$ 
 The above estimate \eqref{zuihou gujishizi} is a direct consequence of Lemma 5.4 in \cite{CZ 21}, Remark~\ref{rem 36} and a well known inequality for the Gaussian distribution: $\abs{ {\Phi}(z_1)- {\Phi}(z_2)}\le \abs{z_1-z_2} $.  {\hfill\large{$\Box$}}\\
\section{Appendix}
In this section, we will make use of the following inequalities repeatedly. We ignore the proof since it is very elementary.
\begin{lemma} \label{upper bound F}
 Assume $\beta>0$ and $\theta>0$. Denote $$A(s)= \int_0^{s} e^{-\theta r} r^{\beta -1}\dif r,\qquad \bar{A}(s)=e^{-\theta s}\int_0^{s} e^{\theta r} r^{\beta -1}\dif r. $$ Then 
there exists a constant $C>0$ such that for any  $s\in [0,\infty)$,
\begin{align*}
A(s)&\le C \times\big(s^{\beta}\mathbbm{1}_{[0,1]}(s) +  \mathbbm{1}_{ (1,\,\infty)}(s)\big)\le C \times(1\wedge s^{\beta}), \\
\bar{A}(s)&\le C \times\big(s^{\beta}\mathbbm{1}_{[0,1]}(s) + s^{\beta-1}\mathbbm{1}_{ (1,\,\infty)}(s)\big)\le  C \times \big (s^{\beta-1} \wedge s^{\beta} \big).
\end{align*}
Especially, when $\beta\in (0,1)$, there exists a constant $C>0$ such that for any  $s\in [0,\infty)$,
\begin{align*}
\bar{A}(s)&\le C \times(1\wedge s^{\beta-1}),\\
\int_0^{\infty} e^{-\theta \abs{r-s}} r^{\beta-1}\dif r&\le C \times (1\wedge s^{\beta-1}),\\
\int_0^{s } e^{-\theta \abs{r-t}} (s-r)^{\beta-1}\dif r&\le C \times (1\wedge (s-t)^{\beta-1}),\qquad t\in [0,s].
\end{align*}
\end{lemma}The following two corollaries are consequences of Lemma~\ref{upper bound F}. We ignore the first one's proof since it is very  direct and simple.
\begin{corollary}\label{ligejifen inequality}
When $H\in (0,\frac12)$, there exists a positive constant $C$ independent on $T$ such that for any $t\in [0,T]$,
\begin{align*}
\int_0^{T } e^{-  \abs{r-t}}  {\frac{\partial   }{\partial r}R^B(r,\,T) }\dif r&\le C \times [1\wedge t^{2H-1} +1\wedge (T-t)^{2H-1}], 
\end{align*}
 and for any $s\in [0,T]$,
\begin{align*}
  \int_0^s e^{r-s}   {\frac{\partial   }{\partial r}R^B(r,\,T)}\dif r &\le C \times[1\wedge s^{2H-1} +1\wedge (T-s)^{2H-1}]\le C\times   {\frac{\partial   }{\partial s}R^B(s,\,T) },\\
 \int_0^s e^{r-T} \abs{\frac{\partial   }{\partial r}R^B(r,\,s)}\dif r &\le C \times[1\wedge s^{2H-1} +1\wedge (T-s)^{2H-1}] \le C\times  {\frac{\partial   }{\partial s}R^B(s,\,T) },\\
 \int_s^T e^{r-T} \abs{\frac{\partial   }{\partial r}R^B(r,\,s) } \dif r&\le C\times  \big[1\wedge (T-s)^{2H-1}\big] ,\\
  \int_0^s e^{r-T} \abs{\frac{\partial   }{\partial s}R^B(r,\,s)}\dif r &\le C\times  \big[1\wedge (T-s)^{2H-1}\big] ,\\
  \int_s^T e^{r-T} \abs{\frac{\partial   }{\partial s}R^B(r,\,s) } \dif r&\le C \times[  s^{2H-1} +1\wedge (T-s)^{2H-1}]\le C\times  {\frac{\partial   }{\partial s}R^B(s,\,T) }.
\end{align*}
\end{corollary}

\begin{nott} For any $v, w\in [0,T]$, define 
\begin{align}
 \psi(w, T)&=\mu(e^{-\abs{\cdot-w} }\mathbbm{1}_{ [0,T]}(\cdot))=\int_{0}^{T} e^{-\abs{u-w} } u^{H-1}\dif u ,\label{psivT defn}\\
 \phi(v, T)&= \int_0^T  \psi(w, T)  \abs{\frac{\partial   }{\partial v}R^B(v,\,w) }  \dif w=\int_{[0,T]^2} e^{-\abs{u-w}}u^{H-1} \abs{\frac{\partial   }{\partial v}R^B(v,\,w) }  \dif u \dif w. \label{phiv1 defn}
\end{align}
\end{nott}
It  is clear that 
 \begin{align}\label{psi vt parital deriva}
 \frac{\partial }{\partial T}\psi(w,T)= e^{w-T} T^{H-1},
\end{align} and  Lemma~\ref{upper bound F} implies that there exists a positive constant $C$ independent on $T$ such that 
{\begin{align}
\psi(T,T)&\le C T^{H-1}, \label{psitt inequality}\end{align} and 
\begin{align} \int_{0 }^{ T} \psi(w,T) \big((T-w)^{2H-1}+w^{2H-1} \big) \dif w&\le C \times \left\{
      \begin{array}{ll}
T^{(3H-1)_{+}}, & \quad \text{if } H\in (0,\frac13)\cup(\frac13,\frac12),\\
\log T, &\quad \text{if } H= \frac13.
 \end{array}
\right.
  \label{psitt inequality 001}
\end{align}} where $a_{+}=\max\set{a,0}$. 
\begin{corollary}\label{cor partialphivt inequality}
When $H\in (0,\frac12)$, there exists a positive constant $C$ independent on $T$ such that for any $v\in [0,T]$,
\begin{align}\label{partialphivt inequality}
0< {\frac{\partial  \phi(v, T)}{\partial T}}\le C\times T^{H-1} \times  {\frac{\partial   }{\partial v}R^B(v,\,T) }.
\end{align}
\end{corollary}
\begin{proof}It is clear that
 \begin{align}\label{partialphivt}
 \frac{\partial  \phi(v, T)}{\partial T}= \abs{\frac{\partial   }{\partial v}R^B(v,\,T) }  \int_0^T e^{-(T-u)} u^{H-1}\dif u+ T^{H-1}\int_0^T e^{-(T-w)}  \abs{\frac{\partial   }{\partial v}R^B(v,\,w) }\dif w.
 \end{align} Lemma~\ref{upper bound F} and Corollary~\ref{ligejifen inequality} imply that when $H\in (0,\frac12)$, 
 \begin{align*}
  \abs{\frac{\partial   }{\partial v}R^B(v,\,T) }   \int_0^T e^{-(T-u)} u^{H-1}\dif u\le C\times T^{H-1} \times  {\frac{\partial   }{\partial v}R^B(v,\,T) };
 \end{align*}and 
 \begin{align*}
 \int_0^T e^{-(T-w)}  \abs{\frac{\partial   }{\partial v}R^B(v,\,w) }\dif w &= \big(\int_0^v+\int_v^T\big) e^{-(T-w)}  \abs{\frac{\partial   }{\partial v}R^B(v,\,w) }\dif w \\
&\le C \times  {\frac{\partial   }{\partial v}R^B(v,\,T) }.
 \end{align*} By plugging these two inequalities into \eqref{partialphivt}, we obtain the desired \eqref{partialphivt inequality}.
\end{proof}
\begin{lemma} \label{phi} 
Let $B(\cdot,\cdot)$ be the Beta function and $\phi(\cdot, T)$ be given as in \eqref{phiv1 defn}. We have that when $H\in (0,\frac12)$,
 \begin{align}\label{lemmma 3.3 limit}
 \lim_{T\to \infty} \frac{ \phi(T, T) }{T^{3H-1}}  =\lim_{T\to \infty}\frac{1}{T^{3H -1}} \int_0^T  \phi(v, T) e^{v-T} \,\dif v=2\big[B(2H,\,H)H-1\big].
 \end{align}In addition, when $\beta> -1$
 ,  there exist a positive constant $C$ independent on $T\ge 1$ such that
 \begin{align}\label{int1to T}
  \int_1^T  \phi(v, T)v^{\beta} \,\dif v \le C \times \left\{
     \begin{array}{ll}
 T^{H }, & \quad \text{when } \beta \in  (-1,-2H ),\\
T^{H }\log T   , & \quad \text{when } \beta=-2H,\\
T^{3H+\beta },& \quad \text{when } \beta>-2H.
  \end{array}
 \right.  
 \end{align}
 Moreover, when $\beta> -2H$, there exist positive constants $c$ and $C$ independent on $T$ such that
 \begin{align}
 c\le\liminf_{T\to \infty}\frac{1}{T^{3H+\beta}} \int_0^T  \phi(v, T)v^{\beta} \,\dif v\le  \limsup_{T\to \infty}\frac{1}{T^{3H+\beta}} \int_0^T  \phi(v, T)v^{\beta} \,\dif v \le C.
 \end{align} 
\end{lemma}
\begin{proof}It is clear that
 \begin{align}
  \frac{ \phi(T, T) }{T^{3H-1}} 
  &= \frac{H}{T^{3H-1}} \Big(\int_{ 0\le u\le w\le T}+\int_{ 0\le w\le u\le T}\Big) e^{-\abs{u-w}}u^{H-1} ((T-w)^{2H-1}-T^{2H-1})\dif u \dif w, \label{lm 3.3 limt}
 \end{align}
and as $T\to \infty$,
  \begin{align*}
  & \frac{1}{T^{3H-1}} \times T^{2H-1} \Big(\int_{ 0\le u\le w\le T}+\int_{ 0\le w\le u\le T}\Big) e^{-\abs{u-w}}u^{H-1}\dif u \dif w\to \frac2H.
 \end{align*}
Denote the other terms as 
  \begin{align*} 
  & \frac{1}{T^{3H-1}} \Big(\int_{ 0\le u\le w\le T}+\int_{ 0\le w\le u\le T}\Big) e^{-\abs{u-w}}u^{H-1} (T-w)^{2H-1} \dif u \dif w:=J_{1}+J_{2}.
 \end{align*}
The change of variables $a=T-w,\,z=a+u$ and the L'H\^{o}pital's rule imply that
 \begin{align*}
\lim_{T\to \infty} J_{1}&=\lim_{T\to \infty} \frac{1}{e^T T^{3H-1} }  \int_0^{T} e^{z} \dif z\int_0^{z} a^{2H-1}(z-a)^{H-1}\dif a
= B(2H,\,H).
 \end{align*}
Similarly, the change of variables $a=T-w$ and the L'H\^{o}pital's rule imply that
\begin{align*}
\lim_{T\to \infty} J_{2}&=\lim_{T\to \infty} \frac{1}{ e^{-T} T^{3H-1}}  \int_0^{T} e^{-u} u^{H-1}\dif u\int_{T-u}^{T} e^{-a} a^{2H-1}\dif a
 =B(2H,\,H).
 \end{align*}
 Substituting the above three limits into \eqref{lm 3.3 limt}, we obtain the first limit of \eqref{lemmma 3.3 limit}.
 
 The L'H\^{o}pital's rule and Coroallry~\ref{cor partialphivt inequality} imply that 
  \begin{align*}
  \lim_{T\to \infty}\frac{1}{T^{3H-1 }} \int_0^T  \phi(v, T)e^{v-T}  \,\dif v &=  \lim_{T\to \infty}\frac{ 1}{e^T T^{3H-1 }} \big[\phi(T, T)e^{T}+  \int_0^T   {\frac{\partial  \phi(v, T)}{\partial T} }e^{v}\,\dif v \big]\\
  &= \lim_{T\to \infty}\frac{   \phi(T, T)}{ T^{3H-1 }}=2\big[B(2H,\,H)H-1\big].
 \end{align*}
 and when $\beta> -2H$,
 \begin{align*}
  \liminf_{T\to \infty}\frac{1}{T^{3H+\beta}} \int_0^T  \phi(v, T)v^{\beta} \,\dif v &\ge c\lim_{T\to \infty}\frac{  \phi(T, T)}{T^{3H-1}}\ge c,\\
 \limsup_{T\to \infty}\frac{1}{T^{3H+\beta}} \int_0^T  \phi(v, T) v^{\beta} \,\dif v &\le C \lim_{T\to \infty}\frac{1}{T^{3H+\beta-1}} \big[{ \phi(T, T)T^{\beta}+ T^{H-1}\int_0^T  \abs{\frac{\partial R^B(v, T)}{\partial v} }v^{\beta}\,\dif v }\big]\\
 &<\infty.
 \end{align*}
With a minor modification of the last limit, we can also obtain the inequality \eqref{int1to T}. 
\end{proof}

\begin{corollary}\label{chen new 3.3}
Let the function $\phi(v, T)$ be given as in \eqref{phiv1 defn}.  
When $H\in (0,\frac12)$, there exists a positive constant $ C$ independent on $T\ge 1$ such that
 \begin{equation}\label{chit shaxiajixian}
\int_0^T \phi(v, T)(v^{2H-1}\wedge 1) \dif v \le C \times \left\{
     \begin{array}{ll}
T^{\gamma_2 }, & \quad \text{when } H \neq \frac14,\\
 T^{\gamma_2 } \log T, & \quad \text{when } H =\frac14,
  \end{array}
 \right.  
 \end{equation}where 
 \begin{equation*}
\gamma_2{=}\left\{
     \begin{array}{ll}
  {H}, & \quad \text{when } H\in (0,\frac14],\\
  {   {5H-1 }}, & \quad \text{when } H\in (\frac14,\frac12).
  \end{array}
 \right.  
\end{equation*}
\end{corollary}
 \begin{proof}
Clearly, 
\begin{align}\label{huafen 2bufen}
 \int_0^T \phi(v, T)(v^{2H-1}\wedge 1) \dif v  = \int_{0}^1  \phi(v,T) \dif v +  \int_{1}^T   \phi(v,T) v^{2H-1}\dif v. \end{align}
 The inequality \eqref{int1to T} implies that when $H\in (0,\frac12)$, 
 \begin{equation}\label{1to T jifen}
    \int_{1}^T   \phi(v,T) v^{2H-1}\dif v\le C  \times \left\{
     \begin{array}{ll}
T^{\gamma_2 }, & \quad \text{when } H \neq \frac14,\\
 T^{\gamma_2 } \log T, & \quad \text{when } H =\frac14,
  \end{array}
 \right.  
\end{equation} 
 By the L'H\^{o}pital's rule and Corollary~\ref{cor partialphivt inequality}, we have for any $T\ge 1$,
  \begin{align*} 
\limsup_{T \to \infty}   \frac{1}{T^{H}}{\int_0^1   \phi(v,T) \dif v } &\le \limsup_{T \to \infty}   \frac{C}{T^{H-1}}{\int_0^1 T^{H-1}\frac{\partial   }{\partial v}R^B(v,\,T)   \dif v }\notag\\
&\le C\limsup_{T \to \infty}   \int_0^1 v^{2H-1}+ (1-v)^{2H-1}  \dif v<\infty ,
  \end{align*} which implies that there exists a positive constant $ C$ independent on $T\ge 1$ such that
\begin{align}\label{2 0to 1 jifen}
\int_0^1   \phi(v,T) \dif v\le C T^{H}.
\end{align}
Plugging the inequalities \eqref{1to T jifen}-\eqref{2 0to 1 jifen} into the identity \eqref{huafen 2bufen}, we obtain the desired estimate \eqref{chit shaxiajixian}.
\end{proof}

\begin{lemma}\label{3.3}
Let the function $\phi(v, T)$ be given as in \eqref{phiv1 defn}.
Denote $$ \chi(T)=\frac{1}{T^{5H-1}}\int_0^T \phi(v, T)  (T-v)^{2H-1}\dif v.$$ 
When $H\in (0,\frac12)$, there exist two positive constant $c,\,C$ independent on $T$ such that
 \begin{align}\label{chit shaxiajixian}
c\le \liminf_{T \to \infty}\chi(T) \le \limsup_{T \to \infty} \chi(T)\le C.
 \end{align}
\end{lemma}
\begin{proof}
First, we rewrite $ \chi(T)$ and then divide the domain $[0,T]^3$ of the integral into several parts by the order of the variables $u,v,w$ as follows:
\begin{align*}
\chi(T)&=\frac{1}{T^{5H-1}}\int_{[0,T]^3} e^{-\abs{u-w}}u^{H-1} \abs{\frac{\partial   }{\partial v}R^B(v,\,w) } (T-v)^{2H-1}\dif u \dif v \dif w\\
 &= \frac{1}{T^{5H-1}}\big(\int_{T>w>v,u>0}+\int_{T>u>w>v>0}+\int_{T>u>v>w>0}+\int_{T>v>w>u>0}+\int_{T>v>u>w>0}\big) \\
 &\times e^{-\abs{u-w}}u^{H-1} \abs{\frac{\partial   }{\partial v}R^B(v,\,w) }(T-v)^{2H-1}\dif v\dif u \dif w\\
 &:=H\times (I_1+I_2+I_3+I_4+I_5),
 \end{align*} where
 \begin{align*}
I_1&=  \frac{1}{T^{5H-1} }\int_{T>w>v,u>0} e^{u-w}u^{H-1} \big(v^{2H-1}+(w-v)^{2H-1}\big)
(T-v)^{2H-1}\dif v\dif u \dif w := I_{11}+I_{12}.
 \end{align*}
By the L'H\^{o}pital's rule, we have
  \begin{align}\label{suplim i11}
\limsup_{T \to \infty}   I_{11}\le \limsup_{T \to \infty}\frac{1}{T^{5H-1}}{\int_0^Te^{-w}\dif w \int^{w}_0 e^{u}u^{H-1}\dif u} {\int^{T}_0v^{2H-1}(T-v)^{2H-1}\dif v} =\frac{B(2H,2H)}H.
  \end{align}
By the change of variables $a=w-v,b=T-v$ and the L'H\^{o}pital's rule again, we have
 \begin{align}
\lim_{T \to \infty} I_{12}&=\lim_{T \to \infty}\frac{1}{e^TT^{5H-1}} {\int_0^Tb^{2H-1}\dif b\int^{b}_0a^{2H-1}\dif a\int^{T-b+a}_0u^{H-1}e^{u+b-a}\dif u}\nonumber\\
&=\lim_{T \to \infty} \frac1{T^{5H-1}}\int_0^Tb^{2H-1}\dif b\int^{b}_0a^{2H-1}(T-b+a)^{H-1}\dif a\notag \\
&= \int_0^1y^{2H-1}\dif y\int^{y}_0 x^{2H-1}(1-y+x)^{H-1}\dif x.
 \end{align}

 Second, in the same vein, we have
 \begin{align}
\limsup_{T \to \infty}  I_2&=\limsup_{T \to \infty} \frac{1}{T^{5H-1}} { \int_{T>u>w>v>0} e^{w-u}u^{H-1} \big(v^{2H-1}+(w-v)^{2H-1}\big)
(T-v)^{2H-1}\dif v\dif u \dif w }\notag \\
 &\le \frac{B(2H,2H)}H+
 \lim_{T \to \infty}  \frac1{e^{-T}T^{5H-1}}{\int_0^Tb^{2H-1}\dif b\int^{b}_0  a^{2H-1}\dif a\int^T_{T+a-b}e^{a-b-u}u^{H-1}\dif u}\nonumber\\
&=\frac{B(2H,2H)}H+ \int_0^1y^{2H-1}\dif y\int^{y}_0 x^{2H-1}(1-y+x)^{H-1}\dif x.
 \end{align}

Third,  it is clear that as $T\to \infty$,
\begin{align}
 I_3& \le \frac{1}{T^{5H-1}} {\int^{T}_0 u^{H-1}\dif u \int_0^{u}e^{v-u}(T-v)^{2H-1}\dif v\int^{v}_0e^{w-v}(v-w)^{2H-1}\dif w}\notag\\
& \le \frac{1}{T^{5H-1}} \int^{T}_0 u^{H-1} (T-u)^{2H-1}\dif u=B(2H,H)\frac{1}{T^{2H}} \to 0,
\end{align} 
and Lemma~\ref{upper bound F} implies that 
\begin{align}
  I_4 &\le  \frac{1}{T^{5H-1}}  \int_0^T  (T-v)^{2H-1}\dif v\int_0^v (v-w)^{2H-1}\dif w \int_0^w e^{u-w}u^{H-1}\dif u \notag\\
 &\le  \frac{C}{T^{5H-1}}  \int_0^T  (T-v)^{2H-1}\dif v\int_0^v (v-w)^{2H-1} w^{H-1}\dif w
=C\times \frac{\Gamma(2H)\Gamma(2H)\Gamma(H)}{\Gamma(5H)},
\end{align}
and 
\begin{align}
  I_5 &\le  \frac{1}{T^{5H-1}}  \int_0^T  (T-v)^{2H-1}\dif v\int_0^v u^{H-1}\dif u  \int_0^u e^{w-u} (v-w)^{2H-1}\dif w\notag\\
 &\le  \frac{1}{T^{5H-1}}  \int_0^T  (T-v)^{2H-1}\dif v\int_0^v (v-u)^{2H-1} u^{H-1}\dif u
=  \frac{\Gamma(2H)\Gamma(2H)\Gamma(H)}{\Gamma(5H)}. \label{suplim i5}
\end{align}
Finally, combining the limits \eqref{suplim i11}-\eqref{suplim i5} together,  we obtain the desired estimate \eqref{chit shaxiajixian}.
\end{proof}

\begin{nott} 
Denote by $\delta_a(\cdot)$ the Dirac delta function centered at a point $a$.  
\end{nott}
\begin{nott}\label{kappa nott}
Let the function $f_T(u,v)$ be given as in \eqref{ft ts 000}.  Denote the $1$-th contraction between $f_T$ and $f_T$ in $\mathfrak{H}_1^{\otimes 2}$ as a function $$\kappa (u_1,u_2) = f_T\otimes_{1'} f_T  .$$\end{nott}

\begin{lemma}\label{lem 3.2 guji ft yasuo}
 When $H\in (0,\frac12)$, there exist a positive constant $C$ independent on $T\ge 1$ such that 
\begin{equation}\label{yasuo bj} 
(\mu\times \mu) (\abs{\kappa}) \le C T^{\gamma_1}, 
\end{equation} 
where for any $\epsilon>0$,
\begin{equation}\label{dingyi gamma1}
\gamma_1=\left\{
      \begin{array}{ll}
 {H}, & \quad \text{when } H\in (0,\frac13),\\
  {H+\epsilon}, & \quad \text{when } H=\frac13,\\
 {4H-1}, &\quad \text{when } H\in (\frac13,\frac12).
 \end{array}
\right.
\end{equation}
\end{lemma}
\begin{proof}
For simplicity, we assume that $\theta=1$.  
First, since 
for any $-\infty<a<b<\infty$,
 \begin{align}\label{distribution deriva}
 \frac{\dif }{\dif x} \mathbf{1}_{[a,b]}(x) =\delta_a(x)-\delta_b(x),
 \end{align}
we have  when $u_1,u_2\in [0,T]$, 
\begin{align*}
 \kappa(u_1,u_2)&=-\int_{[0,T]^2} f_T(u_1,v_1) \frac{\partial}{\partial v_2} f_T(u_2,v_2) \frac{\partial R^B(v_1,v_2)}{\partial v_1} \dif v_1\dif v_2\\
&=-\int_{[0,T]^2} e^{-\abs{u_1-v_1}-\abs{u_2-v_2}} \mathrm{sgn} (u_2-v_2) \frac{\partial R^B(v_1,v_2)}{\partial v_1} \dif v_1\dif v_2\\
&-\int_{[0,T]^2} e^{-\abs{u_1-v_1} - \abs{v_2-u_2}}  \frac{\partial R^B(v_1,v_2)}{\partial v_1} (\delta_0(v_2)-\delta_T(v_2)) \dif v_1\dif v_2\\
&:=H\times (I_1+I_2).
\end{align*}
It is clear that
\begin{align*}
 \abs{I_1}&\le {\int_{0\le v_1\le v_2 \le T} e^{-\abs{u_1-v_1}-\abs{u_2-v_2}}  \big((v_2-v_1)^{2H-1}+v_1^{2H-1} \big) \dif v_1\dif v_2}\\
&+   {\int_{0\le v_2\le v_1\le T} e^{-\abs{u_1-v_1}-\abs{u_2-v_2}}  (v_1-v_2)^{2H-1}  \dif v_1\dif v_2}\\
&:=I_{11}+I_{12},
\end{align*}and
\begin{align*}
 I_2&=e^{- (T-u_2)} \int_{[0,T]} e^{-\abs{u_1-v_1}} \big((T-v_1)^{2H-1}+v_1^{2H-1} \big)  \dif v_1.
\end{align*}
Hence,  we have
\begin{align}\label{mumu yasuo}
(\mu\times \mu)( \abs{\kappa})\le H\times (\mu\times \mu)(I_{11}+I_{12}+I_2 ).
\end{align}

Second,  let $\psi(\cdot,T)$ be given as in \eqref{psivT defn}. 
The L'H\^{o}pital's rule, the identity \eqref{psi vt parital deriva}, the inequalities \eqref{psitt inequality}-\eqref{psitt inequality 001} and Corollary~\ref{ligejifen inequality} imply that when $H\in (0,\frac12)$,
\begin{align}
& \limsup_{T\to \infty}\frac{1}{T^{\gamma_1}} (\mu\times \mu)( I_{11} )\notag \\
 &=\limsup_{T\to \infty}\frac{1}{T^{\gamma_1}}  \int_{0\le v_1\le v_2\le T}  \psi(v_1,T)\psi(v_2,T) \big((v_2-v_1)^{2H-1}+v_1^{2H-1} \big) \dif v_1\dif v_2\label{I11 defn} \\
 &=\limsup_{T\to \infty}\frac{1}{\gamma_1T^{\gamma_1-1}}  \psi(T,T) \int_{0 }^{ T} \psi(v_1,T) \big((T-v_1)^{2H-1}+v_1^{2H-1} \big) \dif v_1 \notag\\
 &+\limsup_{T\to \infty}\frac{1}{\gamma_1  T^{\gamma_1-H}}  \int_{0\le v_1\le v_2\le T} [e^{v_2-T}\psi(v_1,T) +e^{v_1-T}\psi(v_2,T)]\big((v_2-v_1)^{2H-1}+v_1^{2H-1} \big) \dif v_1 \dif v_2 \notag\\
 &\le \limsup_{T\to \infty}\frac{C}{ T^{\gamma_1-H}}  \int_{0 }^{ T} \psi(v,T) \big((T-v)^{2H-1}+v^{2H-1} \big) \dif v \notag\\
 &<\infty, \notag 
\end{align}
and 
\begin{align}
 & \limsup_{T\to \infty} \frac{1}{T^{\gamma_1}} {(\mu\times \mu)( I_{12} ) } \notag \\
&= \limsup_{T\to \infty}  \frac{1}{T^{\gamma_1}}  \int_{0\le v_2\le v_1\le T}  \psi(v_1,T)\psi(v_2,T) \abs{\frac{\partial R^B(v_1,v_2)}{\partial v_1}} \dif v_1\dif v_2 \label{I12 limmt}\\
 &\le \limsup_{T\to \infty}\frac{1}{ \gamma_1T^{\gamma_1-1}}  \psi(T,T) \int_{0 }^{ T} \psi(v_2,T) (T-v_2)^{2H-1} \dif v_2\notag\\
 &+\limsup_{T\to \infty}  \frac{1}{ \gamma_1 T^{\gamma_1-H}}  \int_{0\le v_2\le v_1\le T}  [e^{v_2-T}\psi(v_1,T) +e^{v_1-T}\psi(v_2,T)]  (v_1-v_2)^{2H-1}  \dif v_1\dif v_2\notag\\
 &\le \limsup_{T\to \infty}\frac{C}{ T^{\gamma_1-H}}   \int_{0 }^{ T} \psi(v,T)  (T-v)^{2H-1}  \dif v<\infty,\notag
\end{align}
and
\begin{align*}
 \lim_{T\to \infty} \frac{1}{T^{\gamma_1}}(\mu\times \mu)( I_2 )=
 \lim_{T\to \infty}  \frac{1}{T^{\gamma_1}} \psi(T,T)  \int_{0 }^{ T} \psi(v,T) \big((T-v)^{2H-1}+v^{2H-1} \big) \dif v= 0.
\end{align*} 
Plugging the above three limits into \eqref{mumu yasuo}, we obtain the desired \eqref{yasuo bj}. 
\end{proof}

\begin{lemma}\label{chennew 000}
Let the function $f_T(u,v)$ be given as in \eqref{ft ts 000}. Suppose that $u_1,u_2\in [0,T]$. Denote
\begin{align*}
\varphi(u_1,u_2)= \int_{[0,T]^2} f_T(u_1,v_1) f_T(u_2,v_2) \frac{\partial ^2  }{\partial v_1 \partial v_2}\big(R(v_1,\,v_2) -R^B(v_1,\,v_2)\big)  \dif v_1\dif v_2.
\end{align*} When $H\in (0,\,\frac12)$, there exists a constant $C_{\theta,H} > 0$ such that
\begin{equation}\label{zhou ineq chennew}
 \norm{\varphi }^2_{\mathfrak{H}^{\otimes 2}}\le  {C_{\theta,H}}{T^{2\gamma_1}},
\end{equation} where $\gamma_1$ is given in \eqref{dingyi gamma1}.
\end{lemma}
\begin{proof}First, without loss of generality, we assume $\theta=1$. The inequality \eqref{inner product 00.ineq}  implies that
\begin{align}\label{phi2 h 2}
\norm{\varphi }^2_{\mathfrak{H}^{\otimes 2}}&\le \norm{\varphi }^2_{\mathfrak{H}_1^{\otimes 2}}+\big(C_H' (\mu\times \mu)(\abs{ \varphi}) \big)^2+2C_H' (\mu\times \mu)( \abs{\varphi\otimes_{1'}\varphi}).
\end{align}
Second, it is clear that
\begin{align*}
\norm{\varphi }^2_{\mathfrak{H}_1^{\otimes 2}}&=\int_{[0,T]^4}\frac{ \partial^2}{\partial u_1 \partial w_2} \big(\varphi(u_1,u_2)  \varphi(w_1,w_2) \big)\frac{\partial}{\partial w_1}R^B(u_1,w_1)\frac{\partial}{\partial u_2}R^B(u_2,w_2) \dif {u_1}\dif{u_2}\dif {w_1}\dif{w_2},\end{align*}
and
\begin{align*}
\frac{ \partial }{\partial u_1 } \varphi(u_1,u_2)
&=\int_{[0,T]^2} e^{-\abs{u_1-v_1}-\abs{u_2-v_2}} \frac{\partial ^2  }{\partial v_1 \partial v_2}\big(R(v_1,\,v_2) -R^B(v_1,\,v_2)\big)  \dif v_1\dif v_2\\
&\times \big[  \mathbbm{1}_{[0,T]^2}(u_1,u_2)\mathrm{sgn}(v_1-u_1 ) +\mathbbm{1}_{[0,T]}(u_2) (\delta_0(u_1)-\delta_T(u_1))\big].
\end{align*}
Under Hypothesis~\ref{hypthe 1}, from the assumed inequality \eqref{cond hyp2}, we have
\begin{align}\label{phi2 h1 2}
\norm{\varphi }^2_{\mathfrak{H}_1^{\otimes 2}}&\le C \times (I_1+I_2+I_3),
\end{align}where
\begin{align*}
I_1&=\int_{[0,T]^4} \big(\int_{[0,T]^2} e^{-\abs{u_1-v_1}-\abs{u_2-v_2}} (v_1v_2)^{H-1} \dif v_1\dif v_2\big) \big(\int_{[0,T]^2} e^{-\abs{w_1-v_1'}-\abs{w_2-v_2'}} (v_1'v_2')^{H-1} \dif v_1'\dif v_2'\big)\\
&\times  \abs{\frac{\partial}{\partial w_1}R^B(u_1,w_1)}\abs{\frac{\partial}{\partial u_2}R^B(u_2,w_2) }\dif {u_1}\dif{u_2}\dif {w_1}\dif{w_2}\\
&=\Big[\int_{[0,T]^2} \psi(u_1,T) \psi(w_1,T)\abs{\frac{\partial}{\partial w_1}R^B(u_1,w_1)}\dif u_1\dif w_1 \Big]^2  \le C T^{2\gamma_1},
\end{align*}where $\psi(u,T)$ is given as in \eqref{psivT defn} and the last inequality is from \eqref{I11 defn}-\eqref{I12 limmt}. 
In the same vein, the inequalities \eqref{psitt inequality}-\eqref{psitt inequality 001}  
 and the identities \eqref{I11 defn}-\eqref{I12 limmt} imply that  
\begin{align*}
I_2&=\int_{[0,T]^3} \big(\int_{[0,T]^2} e^{-\abs{u_1-v_1}-\abs{u_2-v_2}} (v_1v_2)^{H-1} \dif v_1\dif v_2\big) \big(\int_{[0,T]^2} e^{-\abs{w_1-v_1'}-(T-v_2')} (v_1'v_2')^{H-1} \dif v_1'\dif v_2'\big)\\
&\times  \abs{\frac{\partial}{\partial w_1}R^B(u_1,w_1)}\abs{\frac{\partial}{\partial u_2}R^B(u_2,T) }\dif {u}_1\dif {u}_2\dif w_1\\
&= \int_{[0,T]^2} \psi(u_1,T) \psi(w_1,T) \abs{\frac{\partial}{\partial w_1}R^B(u_1,w_1)}\dif u_1\dif w_1   \int_{[0,T]}  \psi(u_2,T) \psi(T,T)\abs{\frac{\partial}{\partial u_2}R^B(u_2,T)}\dif u_2\\
&\le C T^{2\gamma_1-1},
\end{align*} and
\begin{align*}
I_3&=\int_{[0,T]^2} \big(\int_{[0,T]^2} e^{- (T-v_1)-\abs{u_2-v_2}} (v_1v_2)^{H-1} \dif v_1\dif v_2\big) \big(\int_{[0,T]^2} e^{-\abs{w_1-v_1'}-(T-v_2')} (v_1'v_2')^{H-1} \dif v_1'\dif v_2'\big)\\
&\times  \abs{\frac{\partial}{\partial w_1}R^B(T,w_1)}\abs{\frac{\partial}{\partial u_2}R^B(u_2,T) }\dif {u}_2\dif w_1\\
&=  \Big[  \int_{[0,T]}  \psi(u_2,T) \psi(T,T)\abs{\frac{\partial}{\partial u_2}R^B(u_2,T)}\dif u_2\Big]^2 \le C T^{2(\gamma_1-1)}.
\end{align*}Plugging the above three estimates into \eqref{phi2 h1 2}, we have
\begin{align}\label{phi2 h1 bound}
\norm{\varphi }^2_{\mathfrak{H}_1^{\otimes 2}}&\le C T^{2\gamma_1}.
\end{align}
Third, it is clear that when $H\in (0,\frac12)$
\begin{align}
(\mu\times \mu)(\abs{ \varphi})&\le C\times \int_{[0,T]^2} \Big( \int_{[0,T]^2} e^{-\abs{u_1-v_1}-\abs{u_2-v_2}} (v_1v_2)^{H-1} \dif v_1\dif v_2\Big)  (u_1 u_2)^{H-1} \dif u_1 \dif u_2 \nonumber\\
&=C\times \Big[\int_{[0,T]^2} e^{-\abs{u_1-v_1} }(u_1v_1)^{H-1} \dif u_1 \dif v_1\Big]^2 \notag\\
&\le C \times \Big[ \int_{0}^T (1\wedge v_1^{H-1}) \,v_1^{H-1} \dif v_1\Big]^2\le C, \label{zhoujianjieguo 000}
\end{align} where we use Lemma~\ref{upper bound F} in the last line.
It is clear that
\begin{align*}
(\varphi\otimes_{1'}\varphi)(u_2,w_2) &=\int_{[0,T]^2}\frac{ \partial }{\partial u_1 }  \phi(u_1,u_2)  \phi(w_1,w_2) \frac{\partial}{\partial w_1}R^B(u_1,w_1) \dif {u}_1\dif {w}_1 , \end{align*}
and
\begin{align*}
(\mu\times \mu)( \abs{\varphi\otimes_{1'}\varphi})&\le C\times (J_1+J_2),
\end{align*}
where
\begin{align*}
J_1&=\int_{[0,T]^4} \big(\int_{[0,T]^2} e^{-\abs{u_1-v_1}-\abs{u_2-v_2}} (v_1v_2)^{H-1} \dif v_1\dif v_2\big) \big(\int_{[0,T]^2} e^{-\abs{w_1-v_1'}-\abs{w_2-v_2'}} (v_1'v_2')^{H-1} \dif v_1'\dif v_2'\big)\\
&\times  \abs{\frac{\partial}{\partial w_1}R^B(u_1,w_1)}  (u_2w_2)^{H-1} \dif {u_1}\dif{u_2}\dif {w_1}\dif{w_2}\\
&\le C\times \int_{[0,T]^2} \psi(u_1,T) \psi(w_1,T) \abs{\frac{\partial}{\partial w_1}R^B(u_1,w_1)}\dif u_1\dif w_1 \le C T^{\gamma_1},
\end{align*} where we use the inequality \eqref{zhoujianjieguo 000} in the last line.  In the same vein,  we have
\begin{align*}
J_2&=\int_{[0,T]^3} \big(\int_{[0,T]^2} e^{- (T-v_1)-\abs{u_2-v_2}} (v_1v_2)^{H-1} \dif v_1\dif v_2\big) \big(\int_{[0,T]^2} e^{-\abs{w_1-v_1'}-\abs{w_2-v_2'}} (v_1'v_2')^{H-1} \dif v_1'\dif v_2'\big)\\
&\times  \abs{\frac{\partial}{\partial w_1}R^B(T,w_1)}  (u_2w_2)^{H-1} \dif {u}_2\dif {w_1}\dif{w_2}\\
&\le C\times  \psi(T,T)\times \int_0^T\psi(w_1,T) \abs{\frac{\partial}{\partial w_1}R^B(T,w_1)} \dif w_1  \le C T^{\gamma_1-1}.
\end{align*}Hence,
\begin{align}\label{mumu phiphi1'}
(\mu\times \mu)( \abs{\varphi\otimes_{1'}\varphi})&\le C\times T^{\gamma_1}.
\end{align}

Finally, substituting \eqref{phi2 h1 bound}-\eqref{mumu phiphi1'} into \eqref{phi2 h 2}, we obtain the desired \eqref{zhou ineq chennew}. 
\end{proof}

 \begin{lemma} \label{lm last zuihou}
 Let $\phi(\cdot, T)$ be given as in \eqref{phiv1 defn}. Denote 
\begin{equation}\label{dingyi gamma}
\gamma=\left\{
       \begin{array}{ll}
  { {4H }}, & \quad \text{when } H\in ( 0,\frac14],\\
  {8H-1}, &\quad \text{when } H\in (\frac14,\frac12).
  \end{array}
 \right.  
\end{equation}
When $H\in (0,\,\frac14)\cup(\frac14,\frac12)$ , there exists a positive constant $C$ independent on $T\ge 1$ such that 
\begin{align}
\int_{T>v_1'>v_1,u,w>0}  e^{-|u-v_1|-\abs{w-v_1'}}  \phi(v_1,T) {\frac{\partial R^B(v_1',\,T)}{\partial v_1'}}  \abs{\frac{\partial   R^B(u,\,w)}{\partial w} } \dif  {v}_1 \dif  {v}_1'\dif w\dif u  &\le C {T^{\gamma-H} } , \label{last lemma 1}\\
 \int_{T>v_1'>v_1,u,w>0}  e^{-|u-v_1|-\abs{w-v_1'}}  \phi(v_1',T) {\frac{\partial R^B(v_1,\,T)}{\partial v_1}}  \abs{\frac{\partial  R^B(u,\,w)  }{\partial w}} \dif  {v}_1 \dif  {v}_1'\dif w\dif u &\le C {T^{\gamma-H} }, \label{last lemma 010}\\
 \int_{T>u>v_1,v_1',w>0}   e^{-\abs{u-v_1}-\abs{w-v_1'}}  \phi(v_1,T) \frac{\partial R^B(v_1',T) }{\partial v_1'}  
\abs{\frac{\partial   }{\partial w}R^B(u,\,w) } \dif  {v}_1\dif v_1' \dif u\dif w &\le C {T^{\gamma-H} },\label{last lemma 011}\\
  \int_{T>u>v_1,v_1',w>0}   e^{-\abs{u-v_1}-\abs{w-v_1'}}\phi(v_1',T) \frac{\partial R^B(v_1,T) }{\partial v_1} 
  \abs{\frac{\partial   }{\partial w}R^B(u,\,w) } \dif  {v}_1\dif v_1' \dif u\dif w &\le C {T^{\gamma-H} }.\label{last lemma 012}
\end{align}
When $H=\frac14$, all of the upper bounds can be replaced by $T^{\gamma-H}\log T$.
\end{lemma}
\begin{proof}
We will  the inequality \eqref{last lemma 1}  firstly. To this end, we claim that there exists a positive constant $C$ independent on $T\ge 1$ such that for any fixed $v_1\in (0,T)$,
\begin{align}\label{guji piaoliang 000}
\int_{T>v_1'>v_1,u,w>0}  e^{-|u-v_1|-\abs{w-v_1'}}  v_1'^{2H-1}   \abs{\frac{\partial   R^B(u,\,w)}{\partial w} }   \dif  {v}_1'\dif w\dif u\le  C   T^{2H}( {v_1}^{2H-1} \wedge 1 )
\end{align} 
and 
\begin{align}\label{guji piaoliang}
  \int_{T>v_1'>v_1,u,w>0}  e^{-|u-v_1|-\abs{w-v_1'}}   (T-v_1')^{2H-1}   \abs{\frac{\partial   R^B(u,\,w)}{\partial w} }   \dif  {v}_1'\dif w\dif u\le C T^{2H}(T-v_1)^{2H-1}. 
\end{align}
 In fact, we can obtain them by dividing $\set{T>v_1'>v_1,u,w>0}$, the domain of the  triple integral, into six sub-domain according to the distinct order of $v_1,u,w$, doing suitable changes of variables and then applying Lemma~\ref{upper bound F}  to these triple integrals.  Since this calculation is very elementary, we ignore the details.  
 
Next,  it follows from the inequalities \eqref{guji piaoliang 000}-\eqref{guji piaoliang}  
that when $H\in (0,\frac12)$,
\begin{align*}
& \int_{T>v_1'>v_1,u,w>0}  e^{-|u-v_1|-\abs{w-v_1'}}  \phi(v_1,T)\abs{\frac{\partial R^B(v_1',\,T)}{\partial v_1'}}  \abs{\frac{\partial   R^B(u,\,w)}{\partial w} } \dif  {v}_1 \dif  {v}_1'\dif w\dif u \nonumber\\
&= H\int_{T>v_1'>v_1,u,w>0}  e^{-|u-v_1|-\abs{w-v_1'}}  \phi(v_1,T) v_1'^{2H-1}   \abs{\frac{\partial   R^B(u,\,w)}{\partial w} } \dif  {v}_1 \dif  {v}_1'\dif w\dif u \nonumber\\
&+H \int_{T>v_1'>v_1,u,w>0}  e^{-|u-v_1|-\abs{w-v_1'}}  \phi(v_1,T) (T-v_1')^{2H-1}   \abs{\frac{\partial   R^B(u,\,w)}{\partial w} } \dif  {v}_1 \dif  {v}_1'\dif w\dif u\nonumber\\
&\le C\times  T^{2H} \int_0^T  \phi(v_1,T) \Big[ ( {v_1}^{2H-1} \wedge 1) +  (T-v_1)^{2H-1} \Big]  \dif v_1  \le C\times {T^{\gamma-H} } , 
\end{align*}  where the last line is from Corollary~\ref{chen new 3.3} and  Lemma~\ref{3.3}.

In the same vein, we have that there exists a positive constant $C$ independent on $T\ge 1$ such that for any fixed $v_1'\in (0,T)$,
\begin{align*} 
e^{-v_1'}\int_{(0, v_1')^3 }  e^{-|u-v_1|+ w }  v_1^{2H-1}   \abs{\frac{\partial   R^B(u,\,w)}{\partial w} }   \dif  {v}_1\dif w\dif u \le C\times   (v_1')^{(4H-1)_{+}} ,
\end{align*} 
where $a_{+}=\max\set{a,0}$, and 
\begin{align*} 
e^{-v_1'}\int_{(0, v_1')^3 }  e^{-|u-v_1|+ w }  (T-v_1)^{2H-1}   \abs{\frac{\partial   R^B(u,\,w)}{\partial w} }   \dif  {v}_1\dif w\dif u\le C \times T^{2H}(T-v_1')^{2H-1},
\end{align*}
which, together with Lemma~\ref{phi} and Lemma~\ref{3.3},  implies that the inequality \eqref{last lemma 010} holds. 

Similarly, the inequalities \eqref{last lemma 011}-\eqref{last lemma 012} are from Corollary~\ref{chen new 3.3} and the following estimates: there exists a positive constant $C$ independent on $T\ge 1$ such that for any fixed  $v_1\in (0,T)$, 
\begin{align*}
& \int_{T>u>v_1,v_1',w>0}   e^{-\abs{u-v_1}-\abs{w-v_1'}}  \frac{\partial R^B(v_1',T) }{\partial v_1'}  
\abs{\frac{\partial   }{\partial w}R^B(u,\,w) }  \dif v_1' \dif u\dif w {\le  C} ,\end{align*} 
and respectively, for any fixed $v_1'\in (0,T)$,
\begin{align*}
& \int_{T>u>v_1,v_1',w>0}   e^{-\abs{u-v_1}-\abs{w-v_1'}}  \frac{\partial R^B(v_1,T) }{\partial v_1} 
  \abs{\frac{\partial   }{\partial w}R^B(u,\,w) } \dif  {v}_1  \dif u\dif w \\ 
  &{\le  C \times T^{2H} \times\Big[(T-v_1')^{2H-1}+ \big( (v_1')^{2H-1}\wedge 1\big) \Big].}
\end{align*}
\end{proof}

\begin{lemma}\label{lemma last 2-kappa}
Let $\kappa (u_1,u_2) $ be given in Notation~\ref{kappa nott}  and  $\gamma$ be given as in \eqref{dingyi gamma}.
There exists a positive constant $C$ independent on $T\ge 1$ such that
\begin{equation}\label{zuinan jiexian shang}
(\mu\times \mu)(\abs{\kappa\otimes_1'\kappa })\le C \times \left\{
       \begin{array}{ll}
 T^{\gamma}, & \quad \text{when } H\neq \frac14,\\
T^{\gamma}\log T, & \quad \text{when } H=\frac14.
  \end{array}
 \right.  
 \end{equation} 
\end{lemma}
\begin{proof} For simplicity, we only show the case of $H\neq \frac14$.
First, recall that 
	\begin{align*}
	\kappa (u_1,u_2) &= - \int_{[0,T]^2} f_T(u_1,v_1) \frac{\partial}{\partial v_2}f_T(u_2,v_2)  \frac{\partial   }{\partial v_1}R^B(v_1,\,v_2) \dif v_1\dif v_2\\
	&=\Big[\int_{[0,T]^2}e^{-\abs{u_1-v_1}-\abs{u_2-v_2}}\mathrm{sgn}(v_2-u_2) \frac{\partial   }{\partial v_1}R^B(v_1,\,v_2)  \dif v_1\dif v_2\\
	&+ \int_{[0,T]}e^{-\abs{u_1-v_1}- (T-u_2)}  \frac{\partial   }{\partial v_1}R^B(v_1,\,T)    \dif v_1\Big] \mathbbm{1}_{[0,T]^2}(u_1,u_2).
	\end{align*}
Hence, 
\begin{align*}
\big(\kappa \otimes_1' \kappa \big)(u_2,w_2 )=-\int_{[0,T]^2}\,\frac{\partial }{\partial u_1}\kappa (u_1,u_2) \kappa (w_1,w_2)\frac{\partial   }{\partial w_1}R^B(u_1,\,w_1)  \dif u_1 \dif w_1 ,
\end{align*}
and similar to the inequality \eqref{mumu yasuo}, we have
\begin{align}\label{fenjiewei 6k}
(\mu\times \mu)(\abs{\kappa \otimes_1' \kappa })\le C \times (K_1+K_2+K_3 +K_4+K_5+K_6),
\end{align}where
\begin{align*}
K_1&=\int_{[0,T]^8}  e^{-\abs{u_1-v_1}-\abs{u_2-v_2}}\abs{ \frac{\partial   }{\partial v_1}R^B(v_1,\,v_2) } e^{-\abs{w_1-v_1'}-\abs{w_2-v_2'}} \abs{\frac{\partial   }{\partial v_1'}R^B(v_1',\,v_2') }  \\
&\times \abs{\frac{\partial   }{\partial w_1}R^B(u_1,\,w_1) } (u_2 w_2)^{H-1}  \dif u_1\dif w_1 \dif u_2\dif w_2  \dif v_1\dif v_2 \dif v_1'\dif v_2' .\\
&=\int_{[0,T]^4} \phi(v_1,T)\phi(v_1',T) e^{-\abs{u-v_1}-\abs{w-v_1'}} \abs{\frac{\partial   }{\partial w}R^B(u,\,w) } \dif  {v}_1\dif v_1' \dif u\dif w.
\end{align*}where $\phi(v_1,T)$ is defined as in \eqref{phiv1 defn}.

We divide the integral $K_1$ by $T^{\gamma}$ and separate the domain of the integration into four parts as follows:
\begin{align}
&\frac1{T^{\gamma}}\int_{[0,T]^4} \phi(v_1,T)\phi(v_1',T) e^{-\abs{u-v_1}-\abs{w-v_1'}} \abs{\frac{\partial   }{\partial w}R^B(u,\,w) } \dif  {v}_1\dif v_1' \dif u\dif w\notag\\
=&\frac1{T^{\gamma}}\big(\int_{T>v_1'>v_1,u,w>0}+\int_{T>v_1>v_1',u,w>0}+\int_{T>u>v_1,v_1',w>0}+\int_{T>w>v_1,u,v_1'>0}\big)\notag\\ 
\times&\phi(v_1,T)\phi(v_1',T) e^{-\abs{u-v_1}-\abs{w-v_1'}} \abs{\frac{\partial   }{\partial w}R^B(u,\,w) } \dif  {v}_1\dif v_1' \dif u\dif w\notag\\
:=& L_1+L_2+L_3+L_4.\label{k1cjuyifenjie}
\end{align}

The L'H\^{o}pital's rule, Lemma~\ref{upper bound F}, Corollary~\ref{ligejifen inequality}, Lemma~\ref{phi} and Corollary~\ref{cor partialphivt inequality} imply that
\begin{align}\label{L1}
& \limsup_{T \to \infty} L_1\\
&\le c\times \limsup_{T \to \infty}  \Big[\frac1{T^{\gamma-3H}e^T}  \int_{[0,T]^3} \phi(v_1,T) e^{-|u-v_1|+w} \abs{\frac{\partial   }{\partial w}R^B(u,\,w) } \dif  {v}_1\dif w\dif u\notag\\
&+\frac1{T^{\gamma-H} } \int_{T>v_1'>v_1,u,w>0}  e^{-|u-v_1|-\abs{w-v_1'}}  \phi(v_1,T)\abs{\frac{\partial R^B(v_1',\,T)}{\partial v_1'}}  \abs{\frac{\partial   }{\partial w}R^B(u,\,w) } \dif  {v}_1 \dif  {v}_1'\dif w\dif u   \notag\\
&+\frac1{T^{\gamma-H} } \int_{T>v_1'>v_1,u,w>0}  e^{-|u-v_1|-\abs{w-v_1'}}  \phi(v_1',T)\abs{\frac{\partial R^B(v_1,\,T)}{\partial v_1}} \abs{\frac{\partial   }{\partial w}R^B(u,\,w) } \dif  {v}_1 \dif  {v}_1'\dif w\dif u  \Big]\notag\\
&:=L_{11}+L_{12}+L_{13}+L_{14}+L_{15},
\end{align} where 
\begin{align}
L_{11}&= \lim_{T\to \infty}\frac1{e^T T^{\gamma-3H}} \int_{T>u>v_1,w>0}\phi(v_1,T)  e^{-|u-v_1|+w} \abs{\frac{\partial   }{\partial w}R^B(u,\,w) } \dif  {v}_1\dif w\dif u\notag\\
&\le H\times \lim_{T \to \infty}\frac{1}{e^{2T} T^{\gamma-3H}}{\int_{[0,T]^2}\phi(v_1,T)e^{v_1+w}(w^{2H-1}+(T-w)^{2H-1})\dif v_1\dif w}\notag\\
&+C \times\lim_{T\to \infty}\frac1{e^T T^{1+\gamma-4H}} \int_{T>u>v_1,w>0}e^{-|u-v_1|+w}  \frac{\partial R^B(v_1,T)  }{\partial v_1}  \abs{\frac{\partial   }{\partial w}R^B(u,\,w) } \dif  {v}_1\dif w\dif u\notag\\
&=0,
\end{align}
and 
\begin{align*}
L_{12}&= \lim_{T\to \infty}\frac1{e^T T^{\gamma-3H}} \int_{T>v_1>u,w>0}\phi(v_1,T)  e^{-|u-v_1|+w} \abs{\frac{\partial   }{\partial w}R^B(u,\,w) } \dif  {v}_1\dif w\dif u \\
&\le C \times  \lim_{T \to \infty}\frac1{e^{2T}T^{1+\gamma-6H}}{\int_{[0,T]^2}e^{u+w}\abs{\frac{\partial   }{\partial w}R^B(u,\,w) } \dif u\dif w}\\
&+C \times \lim_{T\to \infty}\frac1{e^T T^{1+\gamma-4H}} \int_{T>v_1>u,w>0}  e^{-|u-v_1|+w}  \frac{\partial R^B(v_1,T)  }{\partial v_1}  \abs{\frac{\partial   }{\partial w}R^B(u,\,w) } \dif  {v}_1\dif w\dif u \\
&=0,
\end{align*}
and 
\begin{align}%
L_{13}=&\limsup_{T\to \infty}\frac1{ e^T T^{\gamma-3H}}  \int_{T>w>v_1,u>0} \phi(v_1,T)  e^{-|u-v_1|+w} \abs{\frac{\partial   }{\partial w}R^B(u,\,w) } \dif  {v}_1\dif w\dif u\notag\\
\le &\limsup_{T \to \infty}\frac{H }{T^{\gamma-3H}}  {\int_{[0,T]^2}\phi(v_1,T)e^{-|u-v_1|}(T-u)^{2H-1} \dif u\dif v_1} \notag\\
+& C \times \lim_{T\to \infty}\frac1{e^T T^{1+\gamma-4H}} \int_{T>w>v_1,u>0}  e^{-|u-v_1|+w}  \frac{\partial R^B(v_1,T)  }{\partial v_1}  \abs{\frac{\partial   }{\partial w}R^B(u,\,w) } \dif  {v}_1\dif w\dif u\notag\\
\le &C\times\limsup_{T \to \infty} \frac{\int_{0}^{ T}\phi(v,T) (T-v)^{2H-1} \dif r}{T^{\gamma-3H} }<\infty,\label{L13}
\end{align}
where the last inequality is from Lemma~\ref{3.3}.  Lemma~\ref{lm last zuihou} implies that 
\begin{align*}  
L_{14}+L_{15}&=\limsup_{T \to \infty}  \frac1{T^{\gamma-H} } \int_{T>v_1'>v_1,u,w>0}  e^{-|u-v_1|-\abs{w-v_1'}} \abs{\frac{\partial   }{\partial w}R^B(u,\,w) }  \\
&\times \Big[\phi(v_1,T) {\frac{\partial R^B(v_1',\,T)}{\partial v_1'}} + \phi(v_1',T) {\frac{\partial R^B(v_1,\,T)}{\partial v_1}} \Big]  \dif  {v}_1 \dif  {v}_1'\dif w\dif u   \notag\\
&<\infty.
\end{align*}
Hence,
\begin{equation} \label{L1 guji limsup}
\limsup_{T \to \infty} L_1=L_{11}+L_{12}+L_{13}+L_{14}+L_{15}<\infty.
\end{equation}
In the same vein, we have  \begin{align}\label{L1L2} \limsup_{T \to \infty} L_2<\infty. \end{align}

The L'H\^{o}pital's rule, Corollary~\ref{cor partialphivt inequality} and Lemma~\ref{phi} imply that 
\begin{align}
 &\limsup_{T \to \infty}L_{3}\notag \\
 \le & \limsup_{T \to \infty} \frac1{\gamma T^{\gamma-1}}
\int_0^T \phi(v_1,T) e^{-(T-v_1)}  \dif  {v}_1\int_{ [0,T]^2} \phi(v_1',T) e^{ -\abs{w-v_1'}} \abs{\frac{\partial   }{\partial w}R^B(T,\,w) }\dif v_1'  \dif w\notag\\
 +& c\times  \limsup_{T \to \infty} \frac1{\gamma T^{\gamma-H}}\Big[\int_{T>u>v_1,v_1',w>0} \big(\phi(v_1,T) \frac{\partial }{\partial v_1'} R^B(v_1',T) +\phi(v_1',T) \frac{\partial }{\partial v_1}R^B(v_1,T) \big)\notag\\
 & e^{-\abs{u-v_1}-\abs{w-v_1'}} \abs{\frac{\partial   }{\partial w}R^B(u,\,w) } \dif  {v}_1\dif v_1' \dif u\dif w \Big]\notag\\
\le & \limsup_{T \to \infty}\frac{1}{T^{\gamma-3H}}{\int_{[0,T]^2}\phi(v_1',T) e^{-\abs{w-v_1'}}  (T-w)^{2H-1} \dif v_1' \dif w}+C< \infty ,\label{L3}
\end{align} where we use \eqref{L13}  and the inequalities \eqref{last lemma 011}-\eqref{last lemma 012} in the last line.
In the same vein, we have
\begin{align}\label{L4}
\limsup_{T \to \infty}L_4
&<\infty.
\end{align}
Substituting   the limits \eqref{L1 guji limsup}-\eqref{L4} into \eqref{k1cjuyifenjie}, we have that when $H\in (0,\frac12)$, there exist a positive constant $C$ such that 
\begin{align*}
\limsup_{T\to\infty} \frac{K_1}{T^{\gamma}} =\limsup_{T\to\infty} [L_1+L_2+L_3+L_4] \le C ,
\end{align*}which implies that there exist a positive constant $C$ such that 
\[K_1\le C T^{\gamma}.\]

Second, it is clear that when $H\in (0, \frac12)$,
\begin{align*}
K_2&=\int_{[0,T]^7}  e^{-\abs{T-v_1}-\abs{u_2-v_2}}\abs{ \frac{\partial   }{\partial v_1}R^B(v_1,\,v_2) } e^{-\abs{w_1-v_1'}-\abs{w_2-v_2'}} \abs{\frac{\partial   }{\partial v_1'}R^B(v_1',\,v_2') }  \\
&\times \abs{\frac{\partial   }{\partial w_1}R^B(T,\,w_1) } (u_2 w_2)^{H-1}   \dif w_1 \dif u_2\dif w_2  \dif v_1\dif v_2 \dif v_1'\dif v_2' \\
&=\int_0^T \phi(v_1,T) e^{-({T-v_1})} \dif  {v}_1 \int_{[0,T]^2} \phi(v_1',T) e^{ -\abs{w-v_1'}} \abs{\frac{\partial   }{\partial w}R^B(T,\,w) }\dif v_1'  \dif w,\\
&\le C {T^{\gamma-1}},
\end{align*} where the last line is from the expression of $L_3$ and its limit \eqref{L3}.  

Third, it is clear that there exists a positive constant $C$ independent of $T$ such that
\begin{align*}
\int_0^T  \abs{\frac{\partial   }{\partial w_1}R^B(u_1,\,w_1) } \dif u_1\le C \times T\times  \frac{\partial   }{\partial w_1}R^B( w_1,\,T) ,
\end{align*}
which, combined with Lemma~\ref{upper bound F} and Corollary~\ref{ligejifen inequality}, implies that when $H\in (0, \frac12)$,
\begin{align*}
K_3&=\int_{[0,T]^7}  e^{-\abs{u_1-v_1}-\abs{T-u_2}}\abs{ \frac{\partial   }{\partial v_1}R^B(v_1,\,T) } e^{-\abs{w_1-v_1'}-\abs{w_2-v_2'}} \abs{\frac{\partial   }{\partial v_1'}R^B(v_1',\,v_2') }  \\
&\times \abs{\frac{\partial   }{\partial w_1}R^B(u_1,\,w_1) } (u_2 w_2)^{H-1}   \dif w_1\dif u_1 \dif u_2\dif w_2  \dif v_1 \dif v_1'\dif v_2' .\\
&\le C T^{H-1}\int_{[0,T] }  \phi(v_1',T)  \dif v_1' \int_{[0,T]^3} e^{ -\abs{w_1-v_1'}-\abs{u_1-v_1}} \abs{ \frac{\partial   }{\partial v_1}R^B(v_1,\,T) }\abs{\frac{\partial   }{\partial w_1}R^B(u_1,\,w_1) } \dif  {v}_1 \dif u_1\dif w_1\\
&\le C T^{H-1}\int_0^T \phi(v_1',T)   \dif v_1'  \int_{[0,T]^2} e^{ -\abs{w_1-v_1'}}  \abs{\frac{\partial   }{\partial w_1}R^B(u_1,\,w_1) } \dif u_1 \dif w_1\\
&\le  C T^{H}\int_0^T \phi(v_1',T)   \dif v_1'  \int_{[0,T]} e^{ -\abs{w_1-v_1'}} \frac{\partial   }{\partial w_1}R^B( w_1,\,T)  \dif w_1\\
&\le  C T^{H}\int_0^T \phi(v_1',T)   \dif v_1' \le CT^{\gamma},
\end{align*}  where the last inequality is from Lemma~\ref{phi}.

Fourth, Lemma~\ref{upper bound F} and Corollary~\ref{ligejifen inequality} imply that when $H\in (0, \frac12)$,
\begin{align*}
K_4&=\int_{[0,T]^6}  e^{- ({T-v_1})-\abs{T-u_2}}\abs{ \frac{\partial   }{\partial v_1}R^B(v_1,\,T) } e^{-\abs{w_1-v_1'}-\abs{w_2-v_2'}} \abs{\frac{\partial   }{\partial v_1'}R^B(v_1',\,v_2') }  \\
&\times \abs{\frac{\partial   }{\partial w_1}R^B(T,\,w_1) } (u_2 w_2)^{H-1}   \dif w_1  \dif u_2\dif w_2  \dif v_1 \dif v_1'\dif v_2' .\\
&\le CT^{H-1}\int_{[0,T]^3}  \phi(v_1',T) e^{-(T-v_1)-\abs{w_1-v_1'}}  { \frac{\partial   }{\partial v_1}R^B(v_1,\,T) } {\frac{\partial   }{\partial w_1}R^B(T,\,w_1) }   \dif  {v}_1\dif v_1'  \dif w_1\\
& \le CT^{H-1} \int_{[0,T]}  \phi(v_1',T)  \dif v_1'   \le C T^{\gamma},
\end{align*}where the last inequality is from Lemma~\ref{phi}.

Fifth,  Lemma~\ref{upper bound F}  and Corollary~\ref{ligejifen inequality} imply that when $H\in (0, \frac12)$,
\begin{align*}
K_5&=\int_{[0,T]^6}  e^{- \abs{u_1-v_1}-\abs{T-u_2}}\abs{ \frac{\partial   }{\partial v_1}R^B(v_1,\,T) } e^{-\abs{w_1-v_1'}-(T-w_2 )} \abs{\frac{\partial   }{\partial v_1'}R^B(v_1',\,T) }  \\
&\times \abs{\frac{\partial   }{\partial w_1}R^B(u_1,\,w_1) } (u_2 w_2)^{H-1}   \dif u_1 \dif w_1  \dif u_2\dif w_2  \dif v_1 \dif v_1'\\
&\le CT^{2(H-1)}\int_{[0,T]^2} \abs{\frac{\partial   }{\partial w_1}R^B(u_1,\,w_1) }  \dif u_1  \dif w_1\\
&\le  C T^{\gamma},
\end{align*}  
and
\begin{align*}
K_6&=\int_{[0,T]^5}  e^{- ({T-v_1})-(T-u_2)}\abs{ \frac{\partial   }{\partial v_1}R^B(v_1,\,T) } e^{-\abs{w_1-v_1'}-(T-w_2 )} \abs{\frac{\partial   }{\partial v_1'}R^B(v_1',\,T) }  \\
&\times \abs{\frac{\partial   }{\partial w_1}R^B(T,\,w_1) } (u_2 w_2)^{H-1}     \dif w_1  \dif u_2\dif w_2  \dif v_1 \dif v_1'\\
&\le C T^{2(H-1)}\int_{[0,T]}    \abs{\frac{\partial   }{\partial w_1}R^B(T,\,w_1) }   \dif w_1\\
&\le C T^{\gamma}.
\end{align*} 
Substituting these upper bounds of $K_i,\,i=1,\dots, 6$ into \eqref{fenjiewei 6k}, we obtain the desired \eqref{zuinan jiexian shang}.

 \end{proof}

\begin{lemma} \label{lem phi guji}
Let $0\le s<t\le T$ and $\phi_1(u,v) $ be given as in \eqref{pjhi1uv}. Then there exists a  constant $C>0$ independent of $T$ such that for all $s,t\ge 0$,
\begin{align}
\norm{ \phi_1}_{\FH_1^{\otimes 2}}^2&\le C \big(\abs{t-s}^{4H}+\abs{t-s}^{4H+1}+\abs{t-s}^{4H+2}\big),\label{phi1 budsh 1}\\
(\mu\times \mu)( \abs{\phi_1\otimes_{1'} \phi_1}) &\le C \big(\abs{t-s}^{4H}+\abs{t-s}^{4H+1}\big).\label{phi1 budsh 11}
\end{align}
\end{lemma}
\begin{proof} For simplicity, we can assume that $\theta=1$. Denote $\vec{u}=(u_1,u_2)$ and $\vec{v}=(v_1,v_2)$.
First, by the identity \eqref{distribution deriva}, we have
\begin{align}
\norm{ \phi_1}_{\FH_1^{\otimes 2}}^2&= \int_{[0,T]^4}\frac{\partial^2 }{\partial u_1 \partial v_2}\big[e^{- \abs{u_1-v_1}}e^{- \abs{u_2-v_2}}\mathbb{1}_{[s,t]^4}(u_1,u_2,v_1,v_2) \big]\frac{\partial}{\partial u_2} R^B(u_1,u_2) \frac{\partial }{\partial v_1}R^B(v_1,v_2)\dif \vec{u}\dif \vec{v} \nonumber\\
=&\int_{[0,T]^4}e^{- \abs{u_1-v_1}}e^{- \abs{u_2-v_2}} \mathbb{1}_{[s,t]^2}( u_2,\,v_1 ) \frac{\partial }{\partial u_2}R^B(u_1,u_2) \frac{\partial }{\partial v_1}R^B(v_1,v_2)\nonumber \\
&\times  \Big[\mathbb{1}_{[s,t]^2}( u_1,\,v_2 ) \mathrm{sgn}(u_1-v_1)\mathrm{sgn}(v_2-u_2)-\mathbb{1}_{[s,t] }( u_1)\mathrm{sgn}(u_1-v_1)(\delta_s(v_2)-\delta_t(v_2))\nonumber \\
&-\mathbb{1}_{[s,t] }( v_2  )\mathrm{sgn}(v_2-u_2)(\delta_s(u_1)-\delta_t(u_1))+(\delta_s(v_2)-\delta_t(v_2))(\delta_s(u_1)-\delta_t(u_1))\Big]\dif \vec{u} \dif \vec{v}\nonumber \\
&:=I_1 +I_2 +I_3 +I_4 .\label{f2 decomp}
\end{align} It is cleat that
\begin{align}
\abs{I_1}&\le \int_{[s,t]^4}\abs{\frac{\partial}{\partial u_2} R^B(u_1,u_2)}\abs{ \frac{\partial}{\partial v_1} R^B(v_1,v_2)} \dif \vec{u} \dif \vec{v}=\Big(\int_{[s,t]^2}\abs{\frac{\partial R^B(u_1,u_2)}{\partial u_2}} \dif u_1 \dif u_2 \Big)^2, \label{i1 guji}
\end{align}
and
\begin{align*}
\int_{s\le u_1\le u_2\le t}\abs{\frac{\partial R^B(u_1,u_2)}{\partial u_2}} \dif u_1 \dif u_2&\le H\int_{s\le u_1\le u_2\le t} (u_2-u_1)^{2H-1} \dif u_1 \dif u_2 =\frac{1}{2(2H+1)}(t-s)^{2H+1}.
\end{align*}By the fact $1-x^\beta \leq (1-x)^\beta$ for any $x \in [0, 1]$ and $\beta\in (0,1)$, it is clear that when $H\in (0,\frac12)$,
\begin{align}
\int_{s\le u_2\le u_1\le t}\abs{\frac{\partial R^B(u_1,u_2)}{\partial u_2}} \dif u_1 \dif u_2&= H\int_{s\le u_2\le u_1\le t} \big((u_1-u_2)^{2H-1}+u_2^{2H-1} \big)\dif u_1 \dif u_2\nonumber\\
&\le  \frac{ 1}{ 2H+1} (t-s)^{2H+1} .\label{guji budsh 0001}
\end{align}Plugging the above two formula into \eqref{i1 guji}, we have
\begin{align*}\abs{I_1}\le C \abs{t-s}^{4H+2}.
\end{align*}
In the same vein, we have
\begin{align*} \abs{I_2} +\abs{I_3}+ \abs{I_4} \le C \big(\abs{t-s}^{4H}+\abs{t-s}^{4H+1}\big).
\end{align*}Hence, we obtain the inequality \eqref{phi1 budsh 1}.

Next, similar to the inequality \eqref{mumu yasuo}, we have
\begin{align}\label{phi1 yasuo fenjie}
(\mu\times \mu)( \abs{\phi_1\otimes_{1'} \phi_1}) &\le H\times (\mu\times \mu)(J_{11}+J_{12}+J_2+J_{3}),
\end{align}where
\begin{align*}
J_{11}&={\int_{s\le v_1\le v_2 \le t} e^{-\abs{u_1-v_1}-\abs{u_2-v_2}}  \big((v_2-v_1)^{2H-1}+v_1^{2H-1} \big) \dif v_1\dif v_2}\mathbb{1}_{[s,t]^2 }(u_1, u_2  )\\
J_{12}&= {\int_{s\le v_2\le v_1\le t} e^{-\abs{u_1-v_1}-\abs{u_2-v_2}} \big((v_1-v_2)^{2H-1}-v_1^{2H-1} \big) \dif v_1\dif v_2}\mathbb{1}_{[s,t]^2 }(u_1, u_2  )\\
J_2&=e^{-(t-u_2)} \int_s^t e^{-\abs{u_1-v_1}}\big((t-v_1)^{2H-1}+v_1^{2H-1} \big) \dif v_1\mathbb{1}_{[s,t]^2 }(u_1, u_2  )\\
J_3&=e^{-(u_2-s)} \int_s^t e^{-\abs{u_1-v_1}}\big((v_1-s)^{2H-1}-v_1^{2H-1} \big) \dif v_1\mathbb{1}_{[s,t]^2 }(u_1, u_2  ).
\end{align*}
It is clear that
\begin{align*}
&\quad (\mu\times \mu)(J_{11})\\
&={\int_{s\le v_1\le v_2 \le t} \big(\int_{[s,t]^2}e^{-\abs{u_1-v_1}-\abs{u_2-v_2}} (u_1u_2)^{H-1} \dif u_1\dif u_2\big) \big((v_2-v_1)^{2H-1}+v_1^{2H-1} \big) \dif v_1\dif v_2}\\
&\le \frac{1}{H^2} (t-s)^{2H}{\int_{s\le v_1\le v_2 \le t}  \big((v_2-v_1)^{2H-1}+v_1^{2H-1} \big) \dif v_1\dif v_2} \le  \frac{ 1}{ H^3(2H+1)} (t-s)^{4H+1},
\end{align*}where the last line is from the inequality \eqref{guji budsh 0001}. In the same vein, we have
\begin{align*}
(\mu\times \mu)(J_{12})&\le \frac{1}{2H^3(2H+1)}(t-s)^{4H+1},\\
(\mu\times \mu)(J_2+J_{3})&\le \frac{2}{H^3} (t-s)^{4H}.
\end{align*}Plugging the above three inequalities into \eqref{phi1 yasuo fenjie}, we obtain the desired \eqref{phi1 budsh 11}.
\end{proof}
\begin{lemma}\label{lem phi2 guji}
Let $0\le s<t\le T$ and $\phi_2(u,v)$ be given as in \eqref{pjhi12 uv}. 
Then there exists a  constant $C>0$ independent of $T$ such that for all $s,t\ge 0$,
\begin{align}
\norm{ \phi_2}_{\FH_1^{\otimes 2}}^2&\le C \big(\abs{t-s}^{2H}+\abs{t-s}^{2H+1} \big),\label{phi2 budsh 1}\\
(\mu\times \mu)( \abs{\phi_2\otimes_{1'} \phi_2}) &\le C \big(\abs{t-s}^{H}+\abs{t-s}^{2H}+\abs{t-s}^{2H+1}\big).\label{phi2 budsh 11}
\end{align}
\end{lemma}
\begin{proof}
We only give the sketch of the proof since it is similar to that of Lemma~\ref{lem phi guji}. For simplicity, we can assume that $\theta=1$.  Denote $\vec{u}=(u_1,u_2)$ and $\vec{v}=(v_1,v_2)$. First, by the identity \eqref{distribution deriva}, we have
\begin{align*}
\norm{ \phi_2}_{\FH_1^{\otimes 2}}^2&\le \norm{e^{-  \abs{u-v}}\mathbbm{1}_{\set{0\le u\le s,\, s\le v\le t}} }_{\FH_1^{\otimes 2}}^2\\
&=\int_{[0,T]^4}\frac{\partial^2 }{\partial u_1 \partial v_2}\big[e^{- \abs{u_1-v_1}}e^{- \abs{u_2-v_2}}\mathbbm{1}_{[0,s]^2}(u_1,u_2)\mathbbm{1}_{[s,t]^2}(v_1,v_2) \big]\frac{\partial R^B(u_1,u_2)}{\partial u_2} \frac{\partial R^B(v_1,v_2)}{\partial v_1}\dif \vec{u}\dif \vec{v} \nonumber\\
&=\int_{[0,T]^4}e^{- \abs{u_1-v_1}}e^{- \abs{u_2-v_2}} \mathbbm{1}_{[0,s]}(u_2) \mathbbm{1}_{[s,t]}(  v_1 ) \frac{\partial R^B(u_1,u_2)}{\partial u_2} \frac{\partial R^B(v_1,v_2)}{\partial v_1}\nonumber \\
&\times  \Big[\mathbbm{1}_{[0,s]}(u_1) \mathbbm{1}_{[s,t]}(  v_2 )   \mathrm{sgn}(u_1-v_1)\mathrm{sgn}(v_2-u_2)-\mathbbm{1}_{[0,s] }( u_1)\mathrm{sgn}(u_1-v_1)(\delta_s(v_2)-\delta_t(v_2))\nonumber \\
&-\mathbb{1}_{[s,t] }( v_2  )\mathrm{sgn}(v_2-u_2)(\delta_0(u_1)-\delta_s(u_1))+(\delta_s(v_2)-\delta_t(v_2))(\delta_0(u_1)-\delta_s(u_1))\Big]\dif \vec{u} \dif \vec{v}\nonumber \\
&:=I_1 +I_2 +I_3 +I_4 .
\end{align*}
Corollary~\ref{ligejifen inequality} implies that there exists a  constant $C>0$ independent of $T$ such that
\begin{align*}
\abs{I_1}&\le \int_{[0,s]^2}e^{-(s-u_1)-(s-u_2)}\abs{\frac{\partial R^B(u_1,u_2)}{\partial u_2}}\dif \vec{u} \times \int_{[s,t]^2}e^{-(v_1-s)-(v_2-s)}\abs{ \frac{\partial R^B(v_1,v_2)}{\partial v_1}} \dif \vec{v} \\
&\le C\times \int_{[s,t]^2} \abs{ \frac{\partial R^B(v_1,v_2)}{\partial v_1}} \dif \vec{v}\le C \abs{t-s}^{2H+1}.
\end{align*}
In the same vein, we have  that there exists a  constant $C>0$ independent of $T$ such that
\begin{align*}
\abs{I_2} +\abs{I_3}+ \abs{I_4} \le C \big(\abs{t-s}^{2H}+\abs{t-s}^{2H+1}\big).
\end{align*} Combining the above two estimates together, we obtain \eqref{phi2 budsh 1}.

Next, we have
\begin{align}\label{phi2 yasuo fenjie}
(\mu\times \mu)( \abs{\phi_2\otimes_{1'} \phi_2}) &\le H\times    (\mu\times \mu)(J_{11}+J_{12}+J_{21}+J_{22} ),
\end{align}where
\begin{align*}
J_{11}&={\int_{s\le v_1\le v_2\le t}  e^{- (v_1 -u_1)-(v_2-u_2 ) }  \big(v_1^{2H-1}+(v_2-v_1)^{2H-1}\big)     \dif v_1\dif v_2}\mathbb{1}_{[0,\,s]^2 }(u_1, u_2  )\\
&+{\int_{s\le v_2\le v_1\le t}  e^{- (v_1 -u_1)-(v_2-u_2 ) }   (v_1-v_2)^{2H-1}   \dif v_1\dif v_2}\mathbb{1}_{[0,\,s]^2 }(u_1, u_2  )\\
&+\int_s^t \Big[ e^{- (v_1 -u_1)-(t-u_2 ) } \big(v_1^{2H-1}+(t-v_1)^{2H-1}\big)  +e^{- (v_1 -u_1)-(s-u_2 ) }   ( v_1-s)^{2H-1}  \Big] \dif v_1\mathbb{1}_{[0,\,s]^2 }(u_1, u_2  )\\
J_{12}&=\int_{s} ^t \dif v_1\int_0^s \dif v_2\,  e^{- (v_1 -u_1)-(u_2-v_2 ) }  (v_1-v_2)^{2H-1}\mathbb{1}_{[0,\,s] }(u_1)\mathbb{1}_{[s,t] }(u_2  )\\
&+\int_{s} ^t  e^{- (v_1 -u_1)-(u_2-s ) }  (v_1-s)^{2H-1}\dif v_1  \mathbb{1}_{[0,\,s] }(u_1)\mathbb{1}_{[s,t] }(u_2  ) \\
J_{21}&\le \int_0^s  \dif v_1\int_{s} ^t \dif v_2\,  e^{ v_1 -u_1+ u_2-v_2  }\big( v_1^{2H-1}+ (v_2-v_1)^{2H-1} \big)\mathbb{1}_{[s,t] }(u_1) \mathbb{1}_{[0,\,s] } (u_2  )\\
&+2 \int_0^{s}  e^{ v_1 -u_1 +u_2-s  } \big( v_1^{2H-1}+ (s-v_1)^{2H-1} \big)\dif v_1 \mathbb{1}_{[s,t] }(u_1) \mathbb{1}_{[0,\,s] } (u_2  ) \\
J_{22}
&={\int_{0\le v_1\le v_2\le s}  e^{- (u_1 -v_1)-(u_2-v_2 ) }  \big(v_1^{2H-1}+(v_2-v_1)^{2H-1}\big)     \dif v_1\dif v_2}\mathbb{1}_{[s,\,t]^2 }(u_1, u_2  )\\
&+{\int_{0\le v_2\le v_1\le s}  e^{- (u_1 -v_1)-(u_2-v_2 ) }   (v_1-v_2)^{2H-1}   \dif v_1\dif v_2}\mathbb{1}_{[s,\,t]^2 }(u_1, u_2  )\\
&+\int_0^s e^{- (u_1 -v_1)-(u_2 -s) } \big(v_1^{2H-1}+(s-v_1)^{2H-1}\big)   \dif v_1\mathbb{1}_{[s,\,t]^2 }(u_1, u_2  ).
\end{align*}
Lemma~\ref{upper bound F} and Corollary~\ref{ligejifen inequality} imply that there exists a  constant $C>0$ independent of $T$ such that
\begin{align*}
 (\mu\times \mu)(J_{11})&\le C\big((t-s)^{2H+1} +(t-s)^{2H}\big),\quad  (\mu\times \mu)J_{12} \le C(t-s)^{2H},\\
(\mu\times \mu) J_{21}&\le C(t-s)^{H},\quad   (\mu\times \mu)(J_{22})\le C (t-s)^{2H} .
\end{align*}This combined with the inequality \eqref{phi2 yasuo fenjie} proves the proposition. 
\end{proof}

\begin{lemma} \label{ft ht neiji bounds}
Denote $a_{+}=\max\set{a,0}$.
Let $f_T,\, h_T$ be given in \eqref{ft ts 000}-(\ref{ht ts}) respectively. There exists a constant $C>0$ independent on $T\ge 1$ such that 
\begin{align}  
\abs{\innp{  f_T,\,  h_T}_{\mathfrak{H}^{\otimes 2}}} \le C  \times \left\{
      \begin{array}{ll}
T^{(4H-1)_{+}}, & \quad \text{if } H\in (0,\frac14)\cup(\frac14,\frac12),\\
\log T, &\quad \text{if } H= \frac14.
 \end{array}
\right. 
  \end{align}
\end{lemma}
\begin{proof}
For simplicity, we assume that $\theta=1$. Similar to the proof of Lemma~\ref{lem phi guji} and Lemma~\ref{lem phi2 guji}, we have 
\begin{align*}
\abs{\innp{  f_T,\,  h_T}_{\mathfrak{H}_1^{\otimes 2}}}\le I_1+I_2+I_3+I_4, 
\end{align*}where 
\begin{align*}
I_1&=\int_{[0,T]^4} e^{-\abs{u_1-v_1}-(T-u_2)-(T-v_2)}  \abs{\frac{\partial   }{\partial u_2}R^B(u_1,\,u_2) } \abs{\frac{\partial   }{\partial v_1}R^B(v_1,\,v_2) }  \dif {u}_1\dif u_2\dif {v_1}\dif v_2,\\
I_2&=\int_{[0,T]^3} e^{- (T-v_1)-(T-u_2)-(T-v_2)}  \abs{\frac{\partial   }{\partial u_2}R^B(T,\,u_2) } \abs{\frac{\partial   }{\partial v_1}R^B(v_1,\,v_2) }  \dif {u}_2\dif {v_1}\dif v_2,\\
I_3&=\int_{[0,T]^3} e^{- \abs{u_1-v_1}-(T-u_2) }  \abs{\frac{\partial   }{\partial u_2}R^B(u_1,\,u_2) } \abs{\frac{\partial   }{\partial v_1}R^B(v_1,\,T) }  \dif {u_1}\dif u_2\dif  {v}_1,\\
I_4&=\int_{[0,T]^2} e^{- (T-v_1)-(T-u_2) }  \abs{\frac{\partial   }{\partial u_2}R^B(T,\,u_2) } \abs{\frac{\partial   }{\partial v_1}R^B(v_1,\,T) }  \dif {u_2}\dif  {v_1}.
\end{align*}
It follows from Corollary~\ref{ligejifen inequality}, Lemma~\ref{upper bound F} and Lemma~\ref{phi} that when $H\neq \frac14$, there exists a positive constant $C$ independent on $T$ such that 
\begin{align*}
I_1&\le \int_{[0,T]^2} e^{-\abs{u_1-v_1} }  \frac{\partial   }{\partial u_1}R^B(u_1, T)  {\frac{\partial   }{\partial v_1}R^B(v_1,T) }\dif u_1\dif v_1 \le C T^{(4H-1)_{+}},\\
I_2&\le \int_{[0,T]^2} e^{v_1- T +v_2 -T } \abs{\frac{\partial   }{\partial v_1}R^B(v_1,\,v_2) }  \dif { v_1}\dif { v_2} \le C  \int_0^T e^{ v_2 -T} {\frac{\partial   }{\partial v_2}R^B(T,\,v_2) }  \dif v_2\le C,\\
I_3&\le \int_{[0,T]^2} e^{-\abs{u_1-v_1} }  \frac{\partial   }{\partial u_1}R^B(u_1, T)  {\frac{\partial   }{\partial v_1}R^B(v_1,T) }\dif u_1\dif v_1 \le C T^{(4H-1)_{+}},\\
I_4&\le C  \big(\int_0^T e^{ v_1 -T} \abs{\frac{\partial   }{\partial v_1}R^B(T,\,v_1) }  \dif v_1\big)^2\le C, \end{align*}
and 
\begin{align*}
(\mu\times \mu) (\abs{f_T})&=\int_{[0,T]^2}e^{-\abs{u-v}}(uv)^{H-1}\dif u\dif v\le C\int_0^T v^{H-1} (1\wedge v^{H-1})\dif v\le C,\\
 (\mu\times \mu) (\abs{h_T})&=\big( \int_0^T e^{v-T} v^{H-1}\dif v \big)^2\le C\\
(\mu\times \mu) (\abs{f_T\otimes_{1'}h_T})&\le \big(\int_0^T e^{v_2-T}v_2^{H-1}\dif v_2\big)^2 \int_0^T e^{u_2-T}   \frac{\partial   }{\partial u_2}R^B(u_2, T)\dif u_2\\
&+\int_0^T\phi(u_2,T)e^{u_2-T}\dif u_2 \int_0^T e^{v_2-T}v_2^{H-1}\dif v_2 \le C T^{4H-2}.
\end{align*} where $\phi(\cdot, T)$ is given as in \eqref{phiv1 defn}.

Comparing three values $0,\, (4H-1)_{+}$, and $4H-2$, we see that the largest one is $(4H-1)_{+}$. Hence, 
\begin{align*}
\abs{\innp{  f_T,\,  h_T}_{\mathfrak{H}_1^{\otimes 2}}} + (C_{H}' )^2(\mu \times \mu)(\abs{f_T}) (\mu \times \mu)(\abs{h_T})  + 2C_{H}'(\mu \times \mu)(\abs{f_T\otimes_{1'}h_T})\le CT^{(4H-1)_{+}}.
\end{align*}
This combined with the inequality \eqref{inner product 00.ineq} proves the proposition. The case of $H=\frac14$ is in the same vein.
\end{proof}

\vskip 0.2cm {\small {\bf  Acknowledgements}:
 Y. Chen is supported by NSFC (No.11961033).
}


\end{document}